\documentclass[a4paper]{article}

\usepackage{latexsym}

\newtheorem{thm}{Theorem}[section]
\newtheorem{pro}[thm]{Proposition}
\newtheorem{lem}[thm]{Lemma}
\newtheorem{cor}[thm]{Corollary}
\newtheorem{example}[thm]{Example}

\def\qed{{$\Box$}\medskip}

\def\implies{\Rightarrow}

\def\e{\varepsilon} \def\f{\varphi}  \def\l{\lambda}
\def\p{\psi} \def\r{\rho} \def\s{\sigma} \def\ta{\theta} \def\t{\tau}
\def\om{\omega}   \def\D{\Delta}

\def\bC{{\mathbf C}}  
\def\bR{{\mathbf R}}   \def\bZ{{\mathbf Z}}

\def\Re{\mathrm{Re\,}} \def\End{\mathrm{End}}

\def\ip#1#2{{(\,#1\, , \, #2\,)}}
\def\m#1{{\vert #1 \vert}} \def\n#1{\Vert #1 \Vert}
\def\demo{\noindent \emph{Proof.\ }}

\def\AA{{(A,\D)}} \def\aqg{algebraic quantum group}
\def\ar{A_{\mathrm{r}}} \def\dar{\Delta_{\mathrm{r}}}
\def\LL{\Lambda} 
\def\Lim{\mathrm{Lim}} \def\dop{\Delta_{\mathrm{op}}}
\def\wop{W_{\mathrm{op}}}

\def\functions{\mathbf{F}} \def\fs{\mathbf{K}}

\def\End{\hbox{\rm End}}

\def\A{{\cal A}}

\def\implies{\Rightarrow}

\begin{document}

\title{Amenability and Co-Amenability of Algebraic Quantum Groups}

\author{E.~B\'edos$^*$\ \   G.J.~Murphy  \ \ L.~Tuset$^*$}

\date{June 2001 - minor revisions in October 2001 and January 2002}
\maketitle

\begin{abstract}
    
We define concepts of amenability and co-amenability for
algebraic quantum groups in the sense of~A.~Van Daele
\cite{VD}. We show that co-amenability of an \aqg\ always
implies amenability of its dual. Various necessary and/or
sufficient conditions for amenability or co-amenability are
obtained. Co-amenability is shown to have interesting
consequences for the modular theory in the case that the \aqg\
is of compact type.

\vspace{5ex}

\hspace{-\parindent} \emph{Subj. Class.}: Ouantum groups, C*-algebras

\hspace{-\parindent} \emph{MSC 2000}:  Primary 46L05, 46L65.
Secondary 16W30, 22D25, 58B32.

\hspace{-\parindent} \emph{Keywords}: quantum group, amenability.

\end{abstract}

\vfill
\thanks{
$^*$  Partially supported by the Norwegian Research Council}

\newpage

\section{Introduction}

The concept of amenability was first introduced into the realm
of quantum groups by D. Voiculescu in \cite{Vo} for Kac
algebras, and studied further by M. Enock and J.-M. Schwartz
in \cite{ES} and by Z.-J. Ruan in \cite{Ruan} (which also
deals with Hopf von Neumann algebras). On the other hand,
amenability and co-amenability of (regular) multiplicative
unitaries has been defined by S. Baaj and G. Skandalis in
\cite{BS}, and studied in \cite{Ba, Ba2, Bl}.  Amenability and
co-amenability for Hopf C$^{*}$-algebras has been considered
by C.-K. Ng in \cite{Ng1, Ng2}.

The present paper is a continuation of the earlier
paper~\cite{BMT} of the authors, in which we studied the
concept of co-amenability for compact quantum groups as
defined by S.L.~Woronowicz~\cite{MT, Wo}. We showed there that
the quantum group $SU_q(2)$ is co-amenable and that a
co-amenable compact quantum group has a faithful Haar
integral. Combining these results gives a new proof of Nagy's
theorem that the Haar integral of $SU_q(2)$ is
faithful~\cite{Na}.

In this paper we extend the class of quantum groups for which
we study the concept of co-amenability, and we also initiate a
study of the ``dual'' notion of amenability. The quantum groups
we consider are the \aqg s introduced by A.~Van~Daele
in~\cite{VD}. This class is sufficiently large to include
compact quantum groups and discrete quantum groups (for a more
precise statement of what is meant here, see
Proposition~\ref{thm: compact and discrete types} below). An
\aqg\ admits a dual that is also an \aqg; moreover, there is a
Pontryagin-type duality theorem to the effect that the double
dual is canonically isomorphic to the original \aqg\ (see
Proposition~\ref{thm: duality}).

If $\Gamma$ is a discrete group, there are associated to it in
a natural way two \aqg s, namely $(A,\D)=(\bC(\Gamma),\D)$,
and its dual $(\hat A,\hat \D)$, where $\bC(\Gamma)$ is the
group algebra of $\Gamma$ and $\D$ is the co-multiplication on
$\bC(\Gamma)$ given by $\D(x)=x\otimes x$, for all $x\in
\Gamma$. Then $\AA$ is co-amenable if, and only if, $(\hat
A,\hat \D)$ is amenable; moreover, each of these conditions is
in turn equivalent to amenability of $\Gamma$ (see
Examples~\ref{group example} and~\ref{group example
continued}).

We first relate co-amenability of an \aqg\ $\AA$ to a property
of the multiplicative unitary $W$ naturally associated to it.
This provides a link between our theory and that of S.~Baaj
and G.~Skandalis~\cite{BS}, although we make no use of their
results. Then we also obtain several other equivalent
formulations of co-amenability (in Theorem~\ref{thm: net
equivalence}) that generalize well-known results in the group
algebra case.

One basic question in the theory is whether co-amenability of
an \aqg\ $\AA$ is always equivalent to amenability of its dual
$(\hat A,\hat \D)$. In fact, we show that co-amenability of
$\AA$ always implies amenability of $(\hat A,\hat \D)$, but it
is conceivable that the converse is not always true. When
$\AA$ is of compact type (that is, the algebra $A$ is unital)
and its Haar state is tracial, 
the converse is known to hold, as may be deduced from a deep result
of Ruan \cite[Theorem 4.5]{Ruan}.

In another direction, we prove below that if $M$ is the
von~Neumann algebra associated to an \aqg~$\AA$, then co-amenability of
$\AA$ implies injectivity of $M$
(Theorem~\ref{thm:injectivity}).  It may also be deduced from
 \cite[Theorem 4.5]{Ruan} that the converse is true
 in the compact tracial case.
Related to this, we give a direct proof in Theorem~\ref{thm:
injectivity implies dual amenability}, that if $\AA$ is of
compact type and its Haar state is tracial, then
injectivity of~$M$ implies amenability of the dual $(\hat
A,\hat \D)$.

In the final section of this paper we investigate the modular
properties of a co-amenable algebraic quantum group of compact
type. The unital Haar functional $\varphi$ of such an \aqg\ $(A, \D
)$ is a KMS-state when extended to the analytic extension $\ar
$ of $A$. We show, in the case that $\AA$ is co-amenable, the
modular group can be given a description in terms of the
multiplicative unitary of~$\AA$.

We shall continue our investigations of the concepts studied in this paper
in a subsequent paper \cite{BMT2}. Extensions of these results to the
general case of locally compact quantum groups will also be considered.

We now give a  brief summary of how the paper is organized.
Section~2 establishes some preliminaries on multiplier
algebras, C*-algebra and von~Neumann algebras. We give a
careful treatment of slice maps in connection with multiplier
C*-algebras (Theorem~\ref{thm: extended slice map}). This is
material that is often assumed in the literature, but does not
appear to have been anywhere explicitly formulated and
established in terms of known results. Section~3 sets out the
background material on algebraic quantum groups that will be
needed in the sequel. The most important results of the paper
are to be found in Section~4, where most of the results
mentioned in the earlier part of this introduction are proved.
The final section, Section~5, discusses some consequences of
co-amenability for the modular properties of an \aqg\ of
compact type.

We end this introductory section by noting some conventions of
terminology that will be used throughout the paper.

Every algebra will be a (not necessarily unital) associative
algebra over the complex field~$\bC$. The identity map on a
set $V$ will be denoted by $\iota_V$, or simply by $\iota$, if
no ambiguity is involved.

If $V$ and $W$ are linear spaces, $V'$ denotes the linear
space of linear functionals on $V$ and  $V \otimes W$ denotes
the linear space tensor product of $V$ and~$W$. The {\em flip
map}~$\chi$ from ${V \otimes W}$ to ${W \otimes V}$ is the
linear map sending ${v\otimes w}$ onto ${w\otimes v}$, for all
$v\in V$ and $w\in W$. If $V$ and $W$ are Hilbert spaces,
${V\otimes W}$ denotes their Hilbert space tensor product; we
denote by $B(V)$ and $B_0(V)$ the C*-algebras of bounded
linear operators and compact operators on~$V$, respectively.
If $v\in V$ and $w\in W$, $\om_{v,w}$ denotes the weakly
continuous  bounded linear functional on $B(V)$ that maps $x$
onto $\ip{x(v)}w$. We set $\om_v=\om_{v,v}$. We will often
also use the notation $\om_v$ to denote a restriction to a
C*-subalgebra of $B(V)$ (the domain of $\om_v$ will be
determined by the context).

If $V$ and $W$ are algebras, ${V \otimes W}$ denotes their
algebra tensor product. If $V$ and $W$ are C*-algebras, then
${V\otimes W}$ will denote their C*-tensor product with
respect to the minimal (spatial) C*-norm.

\section{C*-algebra preliminaries} \label{sec: preliminaries}

\label{sec: exposition} We shall review in this section  some
results related to multiplier algebras, especially multiplier
algebras of C*-algebras,  and we shall also review elements of
the theory of  (multiplier) slice maps. These topics are
fundamental to C*-algebraic quantum group theory but those
parts of their theory that are most relevant  are scattered
throughout  the literature and are often presented only in a
very sketchy form. Therefore, for the convenience of the
reader and in order to establish notation and terminology, we
present a brief, but sufficiently detailed, background
account.

First, we introduce the multiplier algebra of a non-degenerate
$*$-algebra. This generalizes the usual idea of the multiplier
algebra of a C*-algebra. Recall that a non-zero algebra $A$ is
{\em non-degenerate} if, for every $a \in A$, $a=0$ if $a
b=0$, for all $b\in A$ or  $b a =0$, for all $b\in A$.
Obviously, all unital algebras are non-degenerate. If $A$ and
$B$ are non-degenerate algebras, so is ${A\otimes B}$.

Denote by $\End(A)$ the unital algebra of linear maps from a
non-degenerate $*$-algebra $A$ to itself. Let $M(A)$ denote
the set of elements $x\in{\End}(A)$ such that there exists an
element $y\in {\End}(A)$ satisfying $x(a)^* b = a^* y(b)$, for
all $a,b \in A$. Then $M(A)$ is a unital subalgebra of
${\End}(A)$. The linear map $y$ associated to a given $x\in
M(A)$ is uniquely determined by non-degeneracy and we denote
it by $x^*$. The unital algebra $M(A)$ becomes a $*$-algebra
when endowed with the involution ${x \mapsto x^*}$.

When $A$ is a $C^*$-algebra, the closed graph theorem implies
that $M(A)$ consists of bounded operators. If we endow $M(A)$
with the operator norm, it becomes a C*-algebra. It is then
the usual multiplier algebra in the sense of C*-algebra theory.

Suppose now that $A$ is simply a non-degenerate $*$-algebra
and that $A$ is a self-adjoint ideal in a $*$-algebra~$B$. For
$b \in B$, define $L_b \in M(A)$ by $L_b(a)=b a$, for all $a
\in A$. Then the map $L\colon B \rightarrow M(A)$, $b \mapsto
L_b$, is a $*$-homomorphism. If $A$ is an {\em essential}
ideal in $B$ in the sense that an element $b$ of~$B$ is
necessarily  equal to zero if $b a=0$, for all $a\in A$, or $a
b=0$, for all $a\in A$, then $L$ is injective. In particular,
$A$ is an essential self-adjoint ideal in itself (by
non-degeneracy) and therefore we have an injective
$*$-homomorphism $L\colon A \rightarrow M(A)$. We identify the
image of $A$ under $L$ with $A$. Then $A$ is an essential
self-adjoint ideal of $M(A)$. Obviously, $M(A) = A$ if, and
only if, $A$ is unital.

It $T$ is an arbitrary non-empty set, denote by
$\functions(T)$ and $\fs(T)$ the non-degenerate $*$-algebras
of all complex-valued functions on $T$ and of all
finitely-supported such functions, respectively, where the
operations are the pointwise-defined ones. Clearly, $\fs(T)$
is an essential ideal in $\functions(T)$, and therefore we
have a canonical injective $*$-homomorphism from
$\functions(T)$ to $M(\fs(T))$; a moment's reflection shows
that this homomorphism is surjective and we therefore can, and
do henceforth, use this to identify $M(\fs(T))$ with
$\functions(T)$.

If $A$ and $B$ are non-degenerate $*$-algebras, then it is
easily verified that  $A\otimes B$ is an essential
self-adjoint ideal in $M(A)\otimes M(B)$. Hence, by the
preceding remarks, there exists a canonical injective
$*$-homomorphism from $M(A)\otimes M(B)$ into $M(A\otimes B)$.
We use this to identify $M(A)\otimes M(B)$ as a unital
$*$-subalgebra of $M(A\otimes B)$. In general, these algebras
are not equal.

If $\pi \colon A \to B$ is a homomorphism, it is said to be
{\em non-degenerate} if the linear span of the set
$\pi(A)B={\{\pi(a)b\mid a\in A,\  b\in B\}}$  and the linear span of the set
$B \pi(A)$ are both  equal to~$B$.
In this case, there exists a unique extension to a
homomorphism ${\overline{\pi} \colon M (A) \rightarrow M (B)}$
(see \cite{VD2}), which is determined by $\overline{\pi}(x) (
\pi(a) b) = \pi (x a) b$, for every $ x \in M(A)$, $a \in A$
and $b \in B$. Note that $\overline{\pi}$ is a
$*$-homomorphism whenever $\pi$ is a $*$-homomorphism. We
shall henceforth use the same symbol $\pi$ to denote the
original map and its extension~$\overline{\pi}$.

If $\pi \colon A \to B$ is a $*$-homomorphism between
C*-algebras, we shall use the term non-degenerate only in its
usual sense in C*-theory. Thus, in this case, $\pi$ is
non-degenerate if the closed linear span of the set
$\pi(A)B={\{\pi(a)b\mid a\in A,\  b\in B\}}$ is equal to~$B$.

If $\om$ is a linear functional on~$A$ and $x \in M(A)$, we
define the linear functionals $x \omega$ and $\omega x$ on $A$
by setting $(x \omega)(a) = \omega(a x)$ and $(\omega x)(a) =
\omega(x a)$, for all $a \in A$.

We say $\omega$ is {\em positive} if $\omega(a^*a) \geq 0$,
for all $a \in A$; if $\omega$ is positive, we say it is {\em
faithful} if, for all $a \in A$, $\omega(a^*a) = 0 \Rightarrow
a = 0$.

Suppose given C*-algebras $A$ and $B$. If $\om\in A^*$, the
linear map defined by the assignment ${a\otimes b\mapsto
\om(a)b}$ extends to a norm-bounded linear map ${\om\otimes
\iota}$ from ${A\otimes B}$ to $B$. We call ${\om\otimes
\iota}$ a {\em slice map}. Obviously, if $\t\in B^*$, we can
define the slice map ${\iota\otimes \t\colon A\otimes B\to A}$
in a similar manner. The next result shows how we can extend
these maps to $M(A\otimes B)$. This result is frequently used
in the literature, usually without explicit explanation of how
${\om\otimes \iota}$ is to  be understood or how it is
constructed. Similar remarks apply to the corollary.

\begin{thm} \label{thm: extended slice map}
Let $A$ and $B$ be C*-algebras and let $\om\in A^*$. Then the
slice map ${\om\otimes \iota\colon A\otimes B\to B}$ admits a
unique extension to a norm-bounded linear map ${\om\,\otimes
\,\iota\colon M(A\otimes B)\to M(B)}$ that is strictly
continuous on the unit ball of $M(A\otimes B)$.

If $x\in M(A\otimes B)$ and $b\in B$, then
\[b (\omega \,\otimes\,\iota )(x) =
(\omega \,\otimes\,\iota)((1 \otimes b)x)\quad{\rm and} \quad
(\omega \,\otimes\,\iota )(x)b = (\omega \,\otimes\,\iota
)(x(1 \otimes b)).\]

\end{thm}

\demo To prove this, one may assume that $\om$ is positive,
since the set of positive elements of $A^*$ linearly spans
$A^*$. In this case, the slice map ${\om\otimes \iota}$ is
completely positive~\cite[p.~4]{Wa} and is easily seen to be
strict in the sense defined by E.C.~Lance
in~\cite[p.~49]{Lan}. Hence, by~\cite[Corollary~5.7]{Lan},
${\om\otimes \iota}$ admits an extension to a norm-bounded
linear map ${\om\,\otimes \,\iota\colon M(A\otimes B)\to
M(B)}$ that is strictly continuous on the unit ball of
$M(A\otimes B)$. Uniqueness of ${\om\,\otimes \,\iota}$ and
the properties in the second paragraph of the statement of the
theorem now follow immediately from the strict continuity
condition.~\qed

Of course, an analogous result holds in the preceding theorem
for an element $\t\in B^*$.

Recall that a norm-bounded linear functional on a C*-algebra
$A$ has a unique extension to a norm-bounded strictly
continuous functional on $M(A)$. We shall usually denote the
original functional and its extension by the same symbol. This
result, which is pointed out in the Appendix of \cite{Ku1},
follows easily from a result of D.C.\ Taylor~\cite{Ta} that
asserts that, for each $\om \in A^{*}$, there exist an element
$ a \in A$ and a functional $\theta \in A^{*}$ such that $\om(b)
= \theta ( a b) $, for all $ b \in A$. We then define $\om$ on
$M(A)$ by setting $\om(x) = \theta ( a x) $, for all $ x \in
M(A)$. If $\om \in A^*$ and $\t \in B^*$, it follows that the
norm-bounded linear functional ${\om\otimes \t\colon A\otimes
B\to \bC}$ admits a unique extension to a strictly continuous
norm-bounded linear functional on $M(A\otimes B)$. In
agreement with our standing convention, we shall denote the
extension by the same symbol ${\om\otimes \t}$. Using these
observations, we have the following immediately.

\begin{cor}
Let $A$ and $B$ be C*-algebras and let $\om\in A^*$ and $\t\in
B^*$. Let $x\in M(A\otimes B)$. Then
\[(\om\otimes \t)(x)=\om((\iota\otimes \t)(x))=
\t((\om\otimes \iota)(x)).\]
\end{cor}

We shall also need to consider slice maps in the context of
von~Neumann algebras. Let $M, N$ be von~Neumann algebras on
Hilbert spaces~$H$ and~$K$, respectively. We denote the
von~Neumann algebra tensor product by ${M\,\bar{\otimes}\, N}$
(this is the weak closure of the C*-tensor product ${M \otimes
N}$ in $B(H\otimes K)$). We denote by $M_*$ the predual of $M$
consisting of the normal elements of $M^*$. Recall that for
any $\omega \in M_*$ and $\t\in N_*$, we can define a unique
functional $\omega\bar{\otimes}\t\in (M\bar{\otimes} N)_*$
such that $\|\omega\bar{\otimes}\t\|=\|\omega \|\|\t\|$ and
$(\omega\bar{\otimes}\t )(x\otimes y)=\omega (x)\t (y)$, for
all $x\in M$ and $y\in N$. If $\om\in M_*$, we show now how
one can define a slice map $\om\bar{\otimes}\,\iota $ from
$M\bar{\otimes} N$ to~$N$.
 For any $x\in
M\bar{\otimes} N$, the assignment $\t\mapsto
(\om\bar{\otimes}\t )(x)$ defines a bounded functional on
$N_*$. Since $N_* ^* =N$, there exists a unique element $z\in
N$ such that $(\om\bar{\otimes}\t )(x)=\t (z)$, for all $\t\in N_*$. 
We define
$(\om\bar{\otimes}\iota )(x)$ to be equal to $z$. Thus,
${(\om\bar{\otimes}\t)(x)}=\t((\om\bar{\otimes}\iota )(x))$,
as one would expect of a slice map. Clearly,
$\n{(\om\bar{\otimes}\iota )(x)}\le \n\om\n x$. The map
$\om\bar{\otimes}\,\iota$ which sends $x\in M\bar{\otimes} N$
to $(\om\bar{\otimes}\iota )(x)\in N$ is obviously linear and
norm-bounded. Finally, it is evident that
$\om\bar{\otimes}\,\iota$ is an extension of the usual slice
map ${\om\otimes \iota\colon M\otimes N\to N}$. In a similar
fashion, for each $\t\in N_*$, one can define a slice map
${\iota\bar{\otimes}\,\t\colon M\bar{\otimes} N\to M}$.

We finish this section on C*-algebra preliminaries by
recalling briefly a useful fact concerning completely positive
maps that will be needed in the sequel. Suppose that $\pi : A
\rightarrow B$ is a completely positive unital linear map
between unital C*-algebras~$A$ and~$B$. If $a\in A$ and
$\pi(a^*)\pi(a)=\pi(a^*a)$, then $\pi(x a)=\pi(x)\pi(a)$, for
all $x\in A$~\cite[p.~5]{Wa}. In particular, if $u$ is a
unitary in~$A$ for which $\pi(u)=1$, it follows easily that
$\pi(u^*x u)=\pi(x)$, for all $x\in A$.

\section{Algebraic quantum groups}

We begin this section by defining a multiplier Hopf
$*$-algebra. References for this section are \cite{KD,VD,VD2}.

Let $A$ be a non-degenerate $*$-algebra and let $\D$ be a
non-degenerate $*$-homomorphism from $A$ into $M(A\otimes A)$.
Suppose that the following conditions hold:

(1) $(\D\otimes \iota)\D=(\iota\otimes \D)\D$;

(2) The linear mappings defined by the assignments
 ${a\otimes b\mapsto \D(a)(b\otimes 1)}$ and
${a\otimes b\mapsto \D(a)(1\otimes b)}$ are bijections from
$A\otimes A$ onto itself.

Then the pair $(A,\D)$ is called a {\em multiplier Hopf
$*$-algebra}.

In Condition~(1), we are regarding both maps as maps into
$M(A\otimes A\otimes A)$, so that their equality makes sense.
It follows from Condition~(2), by taking adjoints, that the
maps defined by the assignments
 ${a\otimes b\mapsto (b\otimes 1)\D(a)}$ and
${a\otimes b\mapsto (1\otimes b)\D(a)}$ are also bijections
from $A\otimes A$ onto itself.

Let $\AA$ be a multiplier Hopf $*$-algebra and let $\omega$ be
a linear functional on $A$ and $a$ an element in $A$. There is
a unique element  $(\omega \otimes \iota) \Delta(a)$ in $M(A)$
for which \[(\omega \otimes \iota) (\Delta(a))b = (\omega
\otimes \iota)(\Delta(a)(1 \otimes b))\] and \[ b(\omega
\otimes \iota) (\Delta(a)) = (\omega \otimes \iota)((1 \otimes
b) \Delta(a)),\] for all $b \in A$. The element
$(\iota \otimes \omega) \Delta(a)$ in $M(A)$ is determined similarly.
Thus, $\om$ induces linear maps maps $(\om \otimes \iota)
\Delta$ and $(\iota \otimes \omega) \Delta$ from $A$ to
$M(A)$.

There exists a unique non-zero $*$-homomorphism $\varepsilon$
from $A$ to $\bf C$ such that, for all $a\in A$,
\[(\varepsilon \otimes \iota)\D(a) =(\iota \otimes
\varepsilon)\Delta(a) = a.\] The map $\varepsilon$ is called
the {\em co-unit} of $(A, \Delta)$. Also, there exists a
unique anti-multiplicative linear isomorphism $S$ on $A$ that
satisfies the conditions
\[m(S \otimes \iota) (\Delta(a)(1 \otimes b)) = \varepsilon
(a)b\] \noindent and
\[m(\iota \otimes S)((b\otimes 1)\Delta (a)) = \varepsilon (a)b,\]
\noindent for all $a,b \in A$. Here ${m\colon  A \otimes A\to
A}$ denotes the linearization  of the multiplication map
${A\times A\to A}$. The map $S$ is called the {\em antipode}
of $(A, \Delta)$. Note that $S(S(a^*)^*)=a$, for all $a\in A$.

Let $\pi_{1}$ and $\pi_{2}$ be
non-degenerate homomorphisms from $A$ into some algebras $B$
and $C$, respectively. Clearly, the homomorphism
${\pi_1\otimes \pi_2:A\otimes A\to B\otimes C}$ is then
non-degenerate. Hence, we may form the product ${\pi_{1}
\pi_{2} \colon A \to M(B \otimes C)}$ defined by $\pi_{1}
\pi_{2} = (\pi_{1} \otimes \pi_{2})\D$,  where
$\pi_{1} \otimes \pi_{2}$ is extended to $M(A \otimes A)$ by
non-degeneracy. Obviously, $\pi_1\pi_2$  is a non-degenerate
homomorphism and it  is $*$-preserving whenever both $\pi_{1}$ and~$\pi_{2}$
are $*$-preserving. This product is easily seen to be
associative, with $\varepsilon $ as a unit.

For later use, we remark that the set of non-zero
multiplicative linear functionals $\om$ on $A$ is a group
under this product, with inverse operation given by $\om^{-1}
= \om S$. To see this, note that multiplicativity implies
${\om\otimes \om}=\om\circ m$. Therefore, if $a,b\in A$, we
have $\om(b)(\om(\om S))(a)= {(\om\otimes \om )((b\otimes
1)(\iota\otimes S)(\D(a))}={ \om\circ m(\iota\otimes
S)((b\otimes 1)\D(a))}=\om(\e(a)b)=\om(b)\e(a)$. Hence,
$(\om(\om S)(a))=\e(a)$, for all $a\in A$, as required.

If $\om\in A'$, we say $\omega$ is {\em left
invariant} if $(\iota \otimes \omega ) \Delta(a) =
\omega(a)1$, for all $a \in A$. Right invariance is defined
similarly. If a non-zero left-invariant linear functional
on~$A$ exists, it is unique, up to multiplication by a
non-zero scalar. Similarly, for a non-zero right-invariant
linear functional. If $\f$ is a left-invariant functional
on~$A$, the functional $\p=\f S$ is right invariant.

If $A$ admits a non-zero, left-invariant, positive linear
functional~$\f$, we call $(A, \Delta)$ an {\em algebraic
quantum group} and we call $\f$ a {\em left Haar integral}
on~$\AA$. Faithfulness of $\f$ is automatic.

Note that although $\p=\f S$ is right invariant, it may not be
positive. On the other hand, it is proved in \cite{KD} that a
non-zero, right-invariant, positive linear functional on
$A$---a {\em right Haar integral}---necessarily exists. As for
a left Haar integral, a right Haar integral is necessarily
faithful.

There is a unique bijective homomorphism $\r\colon A\to A$
such that $\f(a b)=\f(b \r(a))$, for all $a,b\in A$. Moreover,
$\r(\r(a^*)^*)=a$.

\vspace{2ex} We now discuss duality of algebraic quantum
groups. If $\AA$ is an \aqg, denote by $\hat A$ the linear
subspace of $A'$ consisting  of all functionals $\varphi a$,
where  $a \in A$. Since $\f a=\r(a)\f$, we have $\hat A={\{a
\varphi\ | \ a \in A \}}$. If $\om_1,\om_2\in \hat A$, one can
define a linear functional ${(\om_1\otimes \om_2)\D}$ on~$A$
by setting ${(\om_1\otimes \om_2)\D(a)}={(\f\otimes
\f)((a_1\otimes a_2)\D(a))}$, where $\om_1=\f a_1$ and
$\om_2=\f a_2$. Using this,  the space $\hat A$ can be made
into a non-degenerate $*$-algebra. The multiplication is given
by $\omega_1 \omega_2 = (\omega_1 \otimes \omega_2)\Delta$ and
the involution is given by setting $\omega^*(a) =
\omega(S(a)^*)^-$, for all $a\in A$ and $\om_1,\om_2,\om\in
\hat A$; it is clear that $\omega_1 \omega_2, \om^*\in A'$ but
one can show that, in fact, $\om_1\om_2, \om^*\in \hat A$.

One can realize $M(\hat A)$ as a linear space by identifying
it as  the linear subspace of $A'$ consisting of all $\om\in
A'$ for which ${(\om\otimes \iota)\D(a)}$ and  ${(\iota\otimes
\om)\D(a)}$ belong to $A$. (It is clear that $\hat A$ belongs
to this subspace.) In this identification of $M(\hat A)$, the
multiplication and involution are determined by
\[ (\om_1\om_2)(a)=\om_1((\iota\otimes
\om_2)\D(a))=\om_2((\om_1\otimes \iota)\D(a))\] and
\[\om^*(a)=\om(S(a)^*)^-,\] for all $a\in A$ and
$\om_1,\om_2,\om\in M(\hat A)$.

Note that the co-unit $\e$ of $A$ is the unit of the
$*$-algebra~$M(\hat A)$.

There is a unique $*$-homomorphism  $\hat \Delta$ from $\hat
A$ to $M(\hat A\otimes \hat A)$ such that for all
$\om_1,\om_2\in \hat A$ and $a,b\in A$,
\[((\om_1\otimes 1)\hat \D(\om_2))(a\otimes b)=(\om_1\otimes
\om_2)(\D(a)(1\otimes b))\] and
\[(\hat \D(\om_1)(1 \otimes \om_2))(a\otimes b)=(\om_1\otimes
\om_2)((a\otimes 1)\D(b)).\] Of course, we are here
identifying ${A'\otimes A'}$ as a linear subspace of
${(A\otimes A)'}$ in the usual way, so that elements of ${\hat
A\otimes \hat A}$ can be regarded as linear functionals on
${A\otimes A}$.

The pair $(\hat A, \hat \Delta)$ is an algebraic quantum
group, called the {\em dual} of $\AA$. Its co-unit $\hat
\varepsilon$ and antipode $\hat S$ are given by $\hat
\varepsilon(a\varphi) = \varphi(a)$ and 
$\hat S(a\varphi)= (a\varphi)\circ S$, for all $a\in A$.

\def\doubledual#1{#1\hat{}\phantom{,}\hat{}}

There is an algebraic quantum group version of Pontryagin's
duality theorem for locally compact abelian groups that
asserts that $\AA$ is canonically isomorphic to the dual of
$(\hat A, \hat \Delta)$; that is, $\AA$ is isomorphic to its
double dual $(\doubledual A,\doubledual \Delta)$. This is
stated more precisely in the following result.

\begin{thm} \label{thm: duality}
Suppose that $\AA$ is an \aqg\ with  double dual $(\doubledual
A,\doubledual \Delta)$. Let ${\pi\colon A\to \doubledual A}$
be the canonical map defined by $\pi(a)(\om)=\om(a)$, for all
$a\in A$ and $\om\in \hat A$. Then $\pi$ is an isomorphism of
the algebraic quantum groups $\AA$ and $(\doubledual A,\doubledual \Delta)$;
that is, $\pi$ is a $*$-algebra isomorphism of $A$
onto~$\doubledual A$ for which $(\pi\otimes
\pi)\Delta=\doubledual\Delta\pi$.
\end{thm}

We shall need to consider an object associated to an \aqg\
called its analytic extension. See \cite{KD} for full details.
We need first to recall the concept of a GNS~pair. Suppose
given a positive linear functional $\omega$ on a
$*$-algebra~$A$. Let $H$ be a Hilbert space, and let
${\LL\colon A\to H}$ be a linear map with dense range for
which $\ip{\LL(a)}{\LL(b)}=\om(b^*a)$, for all $a,b\in A$.
Then we call $(H,\LL)$ a {\em GNS~pair} associated to~$\om$.
As is well known, such a pair always exists and is essentially
unique. For, if $(H',\LL')$ is another GNS~pair associated to
$\om$, the map, ${\LL(a)\mapsto \LL'(a)}$, extends to a
unitary $U:H\to H'$.

If $\f$ is a left Haar integral on an \aqg\ $\AA$, and
$(H,\LL)$ is an associated GNS~pair, then it can be shown that
there is a unique $*$-homomorphism ${\pi\colon A\to B(H)}$
such that $\pi(a)\LL(b)=\LL(a b)$, for all ${a,b\in A}$.
Moreover, $\pi$ is faithful and non-degenerate. We let $\ar $
denote the norm closure of $\pi(A)$ in $B(H)$. Thus, $\ar $ is
a non-degenerate C*-subalgebra of~$B(H)$. The representation
${\pi\colon A\to B(H)}$  is essentially unique, for if
$(H',\LL')$ is another GNS~pair associated to $\f$, and
${\pi'\colon A\to B(H')}$ is the corresponding representation,
then, as we observed above, there exists a unitary ${U\colon
H\to H'}$ such that $U\LL(a)=\LL'(a)$,  for all $a\in A$,  and
consequently, $\pi'(a)=U\pi(a)U^*$.

Now observe that there exists a unique non-degenerate
$*$-homomorphism ${\dar  : \ar  \rightarrow M(\ar  \otimes \ar
)}$ such that, for all $a \in A$ and  all $x \in A \otimes A$,
we have \[{\dar (\pi(a)) (\pi \otimes \pi) (x)} = {(\pi
\otimes \pi) (\Delta (a)x)}\] and \[{(\pi \otimes \pi) (x)\dar
(\pi(a))} = {(\pi \otimes \pi) (x \Delta (a))}.\]

We observe also that if $\om \in \ar^{*}$ and $ x \in \ar $,
then the elements $(\om \otimes \iota)( \dar (x) )$ and $(
\iota \otimes \om)( \dar (x))$ both belong to~$\ar$.

First suppose that $\om$ is given as $\om =
\t(\pi(a)\cdot) $, for some element $a \in A$ and
functional $\t \in \ar^{*}$. For $x = \pi(b)$, where $b \in
A$, we have $ {(\om \otimes \iota)( \dar (x) )} = \pi ( (\t\pi
\otimes \iota )( ( a \otimes 1) \D(b) )) \in \pi (A)$. By
continuity, we get $(\om \otimes \iota)( \dar (x) ) \in \ar $
for all $ x \in \ar$. It now follows that $(\om \otimes
\iota)( \dar (x) )\in \ar$, for arbitrary $\om\in \ar^*$ and
$x\in \ar$, by Taylor's result on linear functionals mentioned
earlier and a continuity argument. That $( \iota \otimes \om)(
\dar (x))\in \ar$ is proved in a similar way.

We also recall that the Banach space $\ar^{*}$ becomes a Banach algebra
under the product induced from $\dar$, i.e. defined by $\t \om = (\t \otimes
\om) \dar$, for all $\t, \om \in \ar^{*}$. 

Since the sets ${\D(A)(1\otimes A)}$ and ${\D(A)(A\otimes 1)}$
span ${A\otimes A}$, ${\dar(\ar)(1\otimes \ar)}$ and
${\dar(\ar)(\ar\otimes 1)}$ have dense linear span in
${\ar\otimes \ar}$. We get from this the following
cancellation laws, for a given functional $\om\in \ar^*$:

1. If $\t\om=0$, for all $\t\in \ar^*$, then $\om=0$;

2. If $\om\t=0$, for all $\t\in \ar^*$, then $\om=0$.

Using these cancellation properties, it follows easily that
\[ \ar = [(\om \otimes \iota)( \dar (x) ) \, | \ x \in \ar, \ \om \in
\ar^{*}] =  [(\iota \otimes \om)( \dar (x) ) \, | \ x \in \ar,
\ \om \in \ar^{*}]. \]
Note that we use $[ \cdot ]$ to denote the closed linear span.

\vspace{1ex}

We also need to recall that there is a unique unitary operator
$W$ on $H \otimes H$ such that
\[W((\Lambda\otimes\Lambda)(\Delta(b)(a \otimes 1)))=
\Lambda (a) \otimes \Lambda (b),\] for all $a, b \in A$. This
unitary satisfies the equation \[W_{12} W_{13} W_{23} =
W_{23} W_{12};\] thus, it is a multiplicative unitary,
said to be {\em associated to}~$(H,\LL)$. 
Here we have used the leg numbering notation of \cite{BS}.

One can show that
$W\in M(\ar \otimes B_0(H))$, so especially $W\in (\ar \otimes B_0(H))^{''}
= M \bar{\otimes} B(H)$, where $M$ denotes the von Neumann algebra 
generated by $\ar$. Further, $\ar $ is the norm
closure of the linear space $\{(\iota \otimes \omega) (W)\mid
{\omega \in B_0(H)^*}\}$. Also,
 $\dar  (a) = {W^* (1 \otimes a) W}$,
for all  $a \in \ar $.

The pair $(\ar, \dar )$ is a reduced locally compact quantum
group in the sense of Definition 4.1 of \cite{KV}; we call it
the {\em analytic extension} of $\AA$ associated to~$\f$.

Consider now the algebraic dual $(\hat{A}, \hat{\D })$ of
$(A,\D )$. A right-invariant linear functional $\hat{\psi}$ is
defined on $\hat{A}$ by setting
$\hat{\psi}(\hat{a})=\varepsilon (a)$, for all $a\in A$. Here
$\hat{a}=a\varphi$ and $\e$ is the co-unit of~$\AA$. Since the
linear map, $A\rightarrow\hat{A}$, $a\mapsto\hat{a}$, is a
bijection (by faithfulness of  $\varphi$), the functional
$\hat{\psi}$ is well defined. Now define a linear map
$\hat{\Lambda} :\hat{A}\rightarrow H$ by setting
$\hat{\Lambda}(\hat {a})=\Lambda (a)$, for all $a\in\A$. Since
$\hat{\psi} (\hat{b}^*\hat{a})=\varphi (b^*
a)=\ip{\LL(a)}{\LL(b)}$, for all $a,b\in A$, it follows that
$(H,\hat{\Lambda} )$ is a GNS-pair associated to $\hat{\psi}$.
It can be shown that it is unitarily equivalent to the
GNS-pair for a  left Haar integral $\hat{\varphi}$ of $(\hat
A,\hat \D)$. Hence, we can use $(H,\hat\LL)$ to define a
representation of the analytic extension $(\hat \ar ,\hat \dar
)$ of $(\hat A, \hat \D)$ on the space~$H$. There is a unique
$*$-homomorphism $\hat{\pi} :\hat{A}\rightarrow B(H)$ such
that  $\hat{\pi} (a )\hat{\Lambda}(b )=\hat{\Lambda}(a b)$,
for all $a,b\in\hat{A}$. Moreover, $\hat\pi$ is faithful and
non-degenerate. Let $\hat A_{\mathrm{r}}$ be the norm closure
of $\hat{\pi} (A)$ in $B(H)$, so $\hat A_{\mathrm{r}}$ is a
non-degenerate C*-subalgebra of $B(H)$. The von Neumann algebra 
generated by $\hat \ar$ will be denoted by $\hat M$. One can show that
$W\in M(B_0(H)\otimes \hat \ar )$ and that $\hat
A_{\mathrm{r}}$ is the norm closure of the linear space ${
\{(\omega\otimes \iota)(W)\mid \omega\in B_0(H )^*\}}$. Define
a linear map $\hat \D_{\mathrm{r}} :\hat
A_{\mathrm{r}}\rightarrow M(\hat A_{\mathrm{r}}\otimes\hat
A_{\mathrm{r}} )$ by setting $\hat \D_{\mathrm{r}} (a )=W
(a\otimes 1)W^*$, for all ${a \in \hat A_{\mathrm{r}}}$. Then
$\hat \D_{\mathrm{r}}$ is  the unique $*$-homomorphism ${\hat
\D_{\mathrm{r}} : \hat A_{\mathrm{r}} \rightarrow M(\hat
A_{\mathrm{r}} \otimes\hat A_{\mathrm{r}})}$ such that, for
all $a \in \hat{A}$ and ${x \in \hat{A} \otimes \hat{A}}$,
\[\hat \D_{\mathrm{r}} (\hat\pi(a)) (\hat\pi \otimes \hat\pi) (x) =
(\hat\pi \otimes \hat\pi) (\hat{\D} (a)x)\] and
\[(\hat\pi \otimes \hat\pi) (x)\hat \D_{\mathrm{r}} (\hat\pi(a))
= (\hat\pi \otimes \hat\pi) (x \hat{\D}(a)).\] Note that one
can show that $W\in M(\ar \otimes \hat \ar )$ and ${(\dar
\otimes \iota)(W)}=W_{13}W_{23}$.

\vspace{2ex} An algebraic quantum group $(A, \Delta)$ is of
{\em compact type} if $A$ is unital, and of {\em discrete
type} if there exists a non-zero element $h \in A$ satisfying
$ah = ha = \varepsilon(a)h$, for all $a \in A$.

\begin{pro} \label{thm: compact and discrete types}
Let $(A, \Delta)$ be an algebraic quantum group. If it is of
compact type, its analytical extension $(\ar , \dar )$ is a
compact quantum group in the sense of Woronowicz. If it is of
discrete type, its analytical extension $(\ar , \dar )$ is a
discrete quantum group in the sense of Woronowicz and Van
Daele.
\end{pro}

The duality of discrete and compact quantum groups is stated
precisely in the following result.

\begin{pro}
An algebraic quantum group $(A, \Delta)$ is of compact type
(respectively, of discrete type) if, and only if, its dual
$(\hat A, \hat \Delta)$ is of discrete type (respectively, of
compact type).
\end{pro}

\begin{example} \label{group example} \rm We finish this
section with a brief discussion of the algebraic quantum
groups associated to a discrete group $\Gamma$. This
illustrates the ideas outlined above and provides the
motivation for concepts we introduce later.

First consider the $*$-algebra $\fs(\Gamma)$. This is provided
with a co-multiplication ${\hat\D}$ making it an algebraic
quantum group by setting $\hat\D(f)(x,y)=f(x y)$, for all
$f\in \fs(\Gamma)$. Here we are identifying
${\fs(\Gamma)\otimes \fs(\Gamma)}$ with $\fs(\Gamma\times
\Gamma)$ by identifying the tensor product ${g\otimes h}$ of
two elements $g,h\in \fs(\Gamma)$ with the function in
$\fs(\Gamma\times \Gamma)$ defined by ${(x,y)\mapsto
g(x)h(y)}$. We then identify ${M(\fs(\Gamma)\otimes
\fs(\Gamma))}$ with ${\functions(\Gamma\times\Gamma)}$. The
reason for using the notation $\hat\D$ will be apparent
shortly.

Now let $A=\bC(\Gamma)$ be the group-algebra of $\Gamma$.
Recall that, as a linear space, $A$ has canonical linear basis
the elements of $\Gamma$ and that the multiplication on $A$
extends that of $\Gamma$ and the adjoint operation is
determined by $x^*=x^{-1}$, for all $x\in \Gamma$. We can make
$A$ into an algebraic quantum group by providing it with the
co-multiplication ${\D:A\to A\otimes A}$ determined on the
elements of $\Gamma$ by setting $\D(x)={x\otimes x}$. We shall
now sketch the proof that the dual $(\hat A,\hat\D)$ is the
algebraic quantum group  $(\fs(\Gamma),\hat\D)$.

First, observe that a left Haar integral for $\AA$ is given by
the unique linear functional $\f$ on $A$ for which
$\f(x)=\delta_{x1}$, for all $x\in \Gamma$, where 1 is the
unit of $\Gamma$ and $\delta$ is the usual Kronecker delta
function. If $x,y\in \Gamma$, then $(x\f)(y)=\f(y
x)=\delta_{x^{-1},y}$. It follows that the functionals $x\f$
($x\in \Gamma$) provide a linear basis for $\hat A$. Hence, if
$e_x$ ($x\in \Gamma$) is the canonical linear basis for
$\fs(\Gamma)$ given by $e_x(y)=\delta_{xy}$, we have a linear
isomorphism from $\hat A$ to $\fs(\Gamma)$ given by mapping
$x\f$ onto $e_{x^{-1}}$. Using this isomorphism as an
identification, it is straightforward to check that the
multiplications, adjoint operations and co-multiplications on
$\hat A$ and $\fs(\Gamma)$ are the same; thus, $(\hat
A,\hat\D)=(\fs(\Gamma),\hat\D)$, as claimed.

A GNS-pair $(H,\Lambda)$ associated to $\f$ is given by taking
$H=\ell^2(\Gamma)$ and $\Lambda(x)=e_{x^{-1}}$, for all $x\in
\Gamma$. We choose $e_{x^{-1}}$ rather than $e_x$ in this
formula so as to ensure $\hat\Lambda(e_x)=e_x$. We need this
to get the correct form for $\hat\pi$: it follows easily now
that the representation, ${\hat\pi:\hat A\to B(H)}$, is the
one obtained by left multiplication by elements of $\hat A$;
it therefore extends from $\hat A=\fs(\Gamma)$ to a
$*$-isomorphism $\hat\pi$ from $\ell^\infty(\Gamma)$ onto a
von~Neumann subalgebra of $B(H)$. It is trivially verified
that this von~Neumann is the one generated by $\hat\pi(\hat
A)$; hence, $\hat\pi(\ell^\infty(\Gamma))=\hat M$.

Of course, the representation $\pi:A\to B(H)$ is the one
associated to the (right) regular representation of $\Gamma$ on
$\ell^2(\Gamma)$. Hence, the analytic extension $\ar$
associated to $\AA$ is the reduced group C*-algebra
$C^*_r(\Gamma)$ and the corresponding von~Neumann algebra $M$
is simply the group von~Neumann algebra of~$\Gamma$.

We shall return to this motivating set-up in the sequel.~\qed
\end{example}

\section{Amenability and co-amenability}

We shall retain all the notation from the preceding section.
If $\AA$ is an \aqg, recall that we use the symbol $M$ to denote the
von~Neumann algebra generated by~$\ar $. Of course, $\ar $ and
$\pi(A)$ are weakly dense in~$M$. Since the map $\dar$ is
unitarily implemented, it has a unique weakly continuous
extension to a unital $*$-homomorphism ${\dar :M\rightarrow
M\bar{\otimes} M}$, given explicitly by $\dar (a )=W^*
(1\otimes a)W$, for all $a\in M $. The Banach space $M_{*}$ may then
be regarded as a Banach
algebra when equipped  with the  canonical multiplication
induced by $\dar$; thus, the product of two elements $\om$ and
$\s$ is given by $\om \s={(\om \bar{\otimes}\s)\circ \dar}$.
 
We use the same symbol $R$ to denote the anti-unitary antipode
of $\ar$ and of~$M$, and we denote by $\tau$ the scaling group
of $(\ar,\dar)$ (see \cite{KD,KV}).

Recall also that we use the symbol $\hat{M}$ to denote the von~Neumann algebra
generated by~$\hat A_{\mathrm{r}}$, so that $\hat
A_{\mathrm{r}}$ and $\hat\pi(A)$ are weakly dense in
$\hat{M}$. As with $\dar$, since $\hat \D_{\mathrm{r}}$ is
unitarily implemented, it has a unique extension to a weakly
continuous unital $*$-homomorphism ${\hat \D_{\mathrm{r}}
:\hat{M}\rightarrow \hat{M}\bar{\otimes}\hat{M}}$, given
explicitly by $\hat\D_{\mathrm{r}}(a)={W(a\otimes 1)W^*}$, for
all $a\in \hat M$.

It should be noted that both $M$ and $\hat{M}$ are in the
standard representation. This follows easily from \cite{KD}
and standard von~Neumann algebra theory (see~\cite{SZ}, for
example). As a consequence the normal states on these algebras
are vector states.

In this section we introduce the concepts of amenability and
co-amenability for an \aqg. We begin with the latter concept.
Our definition is an adaptation of one we gave in~\cite{BMT}
for a compact quantum group.  Suppose then $\AA$ is an \aqg\
and let $(H,\Lambda)$ be a GNS~pair associated to a left Haar
integral. Since the representation ${\pi\colon A\to B(H)}$ is
injective, we can use it  to endow $A$ with a C*-norm by
setting $\n a=\n{\pi(a)}$, for $a\in A$. We say that $(A, \D
)$ is {\em co-amenable} if its co-unit $\e$ is norm-bounded
with respect to this norm.

It follows readily from the remarks in the introduction of our
paper~\cite{BMT} that the  group algebra of a
discrete group $\Gamma$ is co-amenable according to this
definition if, and only if, $\Gamma$ is amenable.

On the other hand, co-amenability is automatic in the case of
a discrete-type \aqg:

\begin{pro}
An algebraic quantum group  of discrete type is co-amenable.
\end{pro}

\demo We may suppose our \aqg\ of discrete type is the dual
$(\hat{A}, \hat{\D })$ of an \aqg\ $\AA$ of compact type. If
$\f$ is a left Haar integral for $\AA$ and $a\in A$, then
$\hat{\varepsilon} (\hat{a})=\varphi (a)$, for all $a\in A$,
where, as usual, $\hat a=a\f$. Now $\n{\hat{\pi} (\hat{a})}
 \geq \n{\hat{\Lambda}
(\hat{a}\hat{1} )}_2$, since $\|\hat{\Lambda}(\hat{1} )\|_2
=\|\Lambda (1)\|_2=\varphi (1^* 1 )^{1/2} =1$. Let $c\in A$
and write  $\D (c)=\sum_i b_i\otimes c_i$. Then, by
definition,
\[(\hat{a}\hat{1})(c)=(a\varphi\odot \varphi )\D (c)=\sum_i
(a\varphi )(b_i )(\varphi )(c_i ) =\sum_i \varphi (b_i
a)\varphi (c_i )\] \[=\varphi (\sum_i b_i\varphi (c_i )a
)=\varphi (\varphi (c)a)=\varphi
(a)\varphi (c)=\f(a)\hat{1 }(c).\] Hence,
$\hat{a}\hat{1}=\f(a)\hat{1 }$ and therefore
$\hat{\Lambda}(\hat{a}\hat{1})=\f(a)\hat{\Lambda}(\hat{1
})=\varphi (a)\Lambda (1)$. This implies
\[\|\hat{\pi} (\hat{a})\|\geq \|\hat{\Lambda}
(\hat{a}\hat{1} )\|_2 =|\varphi (a)|\|\Lambda (1)\|_2
=|\varphi (a)|=|\hat{\varepsilon}(\hat{a})|.\] Since the map,
$a\mapsto\hat{a}$, is a bijection from $A$ onto $\hat{A}$, it
follows that the co-unit $\hat\e$ is norm-decreasing and
therefore $(\hat A,\hat\D)$ is co-amenable.~\qed

Co-amenability may be characterized in several ways. All of
the following are essentially well-known in the case of the
group algebra of a discrete group.  Most of these
characterizations are related to results obtained by various
authors in different settings (see \cite{BMT, Bl, ES, Ng1,
Ng2, Ruan, Vo}).

\begin{thm} \label{thm: net equivalence}
Let $(H,\LL)$ be a GNS~pair, and $W$ the corresponding
multiplicative unitary, associated to an algebraic quantum
group $(A, \D )$. Then the following are equivalent
conditions:

(1) $(A, \D )$ is co-amenable;

(2) There exists a net $(v_i)$ of unit vectors in $H$ such
that \[\lim_i\| W(v_i\otimes v)-v_i\otimes v\|_2 =0,\]
\hspace{9ex}for all $v\in H$;

(3) There exists a state $\e_{\mathrm{r}}$ on $\ar$ such that $(\e_{\mathrm{r}}\otimes
\iota)(W)=1$;

(4) There exists a non-zero multiplicative linear functional
on $\ar$;

(5)  The Banach algebra $\ar^{*}$ is unital;

(6)  The Banach  algebra $M_{*}$ has a bounded  left
approximate unit;

(7) $M_{*}$ has a bounded right approximate unit;

(8) $M_{*}$ has a bounded two-sided approximate unit.

\end{thm}

\demo Note first that Condition~(3) makes sense---that is,
$(\e_{\mathrm{r}}\otimes \iota)(W)$ is defined---since ${W\in M(\ar \otimes
B_0(H))}$.

Now suppose Condition~(1) holds and we shall show
that~(3) follows. Since $\e$ is norm-bounded, there is clearly
a unique norm-bounded multiplicative linear functional
$\e_{\mathrm{r}}$ on $\ar $ such that $\e_{\mathrm{r}}\circ
\pi=\e$. Obviously, ${(\e_{\mathrm{r}}\otimes \iota)\dar
=\iota}.$ Using the argument of \cite[Theorem~2.5]{BMT}, we
show now that ${(\e_{\mathrm{r}}\otimes
\iota)(W)=1}.$ We have 
$$W={((\e_{\mathrm{r}}\otimes \iota)\dar
\otimes \iota)(W)}={(\e_{\mathrm{r}}\otimes \iota\otimes
\iota)(\dar \otimes \iota)(W)}$$
$$={(\e_{\mathrm{r}}\otimes
\iota\otimes
\iota)(W_{13}W_{23})}={(1\otimes (\e_{\mathrm{r}}\otimes
\iota)(W))W}.$$ As $W$ is invertible it follows that
${(\e_{\mathrm{r}}\otimes \iota)(W)=1}$, hence 
Condition~(3) holds.

Conversely, suppose Condition~(3) holds and we shall show
that~(1) follows. Let $\e_{\mathrm{r}}$ be as in (3). Since $\e_{\mathrm{r}}$ is a positive
linear functional on~$\ar,$ the slice map ${\e_{\mathrm{r}}\otimes \iota}$ from 
$M(\ar \otimes B_0(H))$  to $M(B_0(H)) = B(H)$ is
completely positive. Note that for $a\in \ar$, $1 \otimes a \in 
M(\ar \otimes B_0(H)),$ so   $\dar
(a)={W^*(1\otimes a)W}$ belongs to $ M(\ar \otimes B_0(H)).$ Hence,
using (3) and our remark at the end of section 2, we get, for all
$a\in \ar,$ 
$${(\e_{\mathrm{r}}\otimes \iota)\dar (a)}={(\e_{\mathrm{r}}\otimes
\iota)(W^*(1\otimes a)W)}$$ 
$$={(\e_{\mathrm{r}}\otimes \iota)(W^*)(\e_{\mathrm{r}}\otimes
\iota)(1\otimes a)(\e_{\mathrm{r}}\otimes \iota)(W)= 1^{*} a \,1 = a}.$$ Now set
$\delta=\e_{\mathrm{r}}\circ \pi$. Then, for all $a,b\in A$,
$$\pi(((\delta\otimes \iota)\D(a))b)={((\e_{\mathrm{r}}\otimes
\iota)(\pi\otimes \pi)\D(a))\pi(b)}$$
$$= {((\e_{\mathrm{r}}\otimes \iota)\dar
(\pi(a)))\pi(b)}=\pi(a)\pi(b)=\pi(a b).$$ Hence,
$((\delta\otimes \iota)\D(a))b=ab$, by injectivity of~$\pi,$
and therefore $(\delta\otimes \iota)\D(a)=a,$ by
non-degeneracy of $A.$ It follows that $\delta=\e.$ Hence, for
$a\in A,$ $\n{\e(a)}= \n{\e_{\mathrm{r}}(\pi(a)}\le \n{\pi(a)}=\n a.$
Therefore, $\e$ is norm-bounded; that is, $\AA$ is
co-amenable. 

To prove the implication (2) $\implies$ (3), suppose
there exists a net of unit vectors $(v_i)$ such that $\lim_i\|
W(v_i\otimes v)-v_i\otimes v\|_2 =0$, for all $v\in H$. By
weak* compactness of the state space of $B(H)$, the net
$(\omega_{v_i})$ of vector states on~$B(H)$ has a state $\e'$
on~$B(H)$ as  an accumulation point. By going to a subnet of
$(v_i)$, if necessary we may suppose that $\e'(x)=\lim_i (xv_i ,v_i),$
for all $x\in B(H)$. Let $\e_{\mathrm{r}}$ denote the 
restriction of $\e'$ to $\ar.$ In the following,  the slice maps $\e_{\mathrm{r}} 
\otimes \iota$ and $\iota \otimes \om_{v}$ for $v \in H$  are defined 
on $M(\ar \otimes B_0(H)).$  Using the assumption, we get 
$$\om_v (\e_{\mathrm{r}}\otimes
\iota)(W) = \e_{\mathrm{r}} (\iota \otimes \om_{v})(W) = \lim_i ((\iota\otimes\omega_v )(W)v_i ,v_i )$$
$$=\lim_i (W(v_i\otimes v),v_i\otimes v)=\lim_i (v_i\otimes
v,v_i\otimes v)= \om_{v}(1)$$ for all $v\in H.$  It follows that
$(\e_{\mathrm{r}}\otimes\iota )(W)=1,$ hence $\e_{\mathrm{r}}$ satisfies $(3).$

Suppose now Condition~(3) holds, so that there exists a state
$\e_{\mathrm{r}}$ on $\ar$ such that $(\e_{\mathrm{r}}\otimes \iota)(W)=1.$ 
Let $\e'$ be a state extension of $\e_{\mathrm{r}}$ to $M.$ Using 
the well known fact that the set of normal states on $M$ is
weak* dense in the set of states on $M$ in combination with the fact 
that every normal state on $M$ is a vector state (as $M$ is in 
standard form), we deduce that there exists a net $(v_i)$ of unit 
vectors in $H$ such that
$\e'(x)=\lim_i (xv_i ,v_i )$, for all $x \in M .$ Then, for
all $v\in H,$
\[\lim_i (W(v_i\otimes v),v_i\otimes v)=
\lim_i ((\iota\otimes\omega_v )(W)v_i ,v_i)=\e'((\iota\otimes\om_v )(W))\]
\[=\e_{\mathrm{r}}((\iota\otimes\om_v )(W))=\omega_v ((\e_{\mathrm{r}}\otimes\iota )(W)) 
=\omega_v (1) =\lim_i (v_i\otimes v,v_i\otimes v).\]  It is now
straightforward to check that $\lim_i\| W(v_i\otimes
v)-v_i\otimes v\|_2 =0$. This proves Condition~(2) holds.

If Condition (1) holds, then the norm-bounded linear
functional $\e_{\mathrm{r}}$ defined on $\ar$ as above is
obviously non-zero and multiplicative, and it is easily seen
to be a unit for $\ar^{*}$. Hence Conditions~(4) and (5)
follow from~(1).

Suppose Condition (4) holds and let $\eta$ be non-zero
multiplicative linear functional on~$\ar$. It is well known
that such a functional is norm bounded. Using this and the
norm-boundedness of the anti-unitary antipode $R$,  it is then clearly enough to show
that $(\eta \otimes \eta R ) \dar ( \pi (a))  = \e(a)$, for all $a
\in A$, in order to show that Condition (1) holds. First we
show that  that any multiplicative linear functional $\omega $
on~$A$ is invariant under~$S^{2}$. As pointed out in section
3, the set of non-zero multiplicative linear functionals on
$A$ has a group structure such that $\omega^{-1} = \omega S$.
Therefore we get $\omega S^{2} = (\omega^{-1})^{-1} = \omega$,
as required. Now set $\om=\eta\pi$. If $a \in A$, we infer from
\cite[Proposition 5.5]{KD} that $\pi(a)$ is an analytic
element of the scaling group $\tau$ on $\ar$ and $\tau_{ni}
(\pi(a)) = \pi(S^{-2n}(a))$, for every integer~$n$. This
implies that
\[\eta \tau_{ni} (\pi(a)) = \omega (S^{-2n}(a))  = \omega (a) =
\eta
\pi(a) .\] By analyticity of the group $\t$, it follows that
$\eta \tau_{i/2} = \eta$ on $\pi(A)$. (This may be seen as follows. 
It is known \cite[Proposition 4.23]{Ku3} that $\t_{t}$ leaves $\pi(A)$ invariant
for each $t \in \bR.$ As $\pi(A)$ is dense in $\ar,$ Corollary 1.22 
of \cite{Ku2} implies 
that $\pi(A)$ is a core for $\t_{z}$ for any $z \in \bC .$ Hence 
$\eta(\t_{n i}(x) )= \eta(x)$ for all $n \in \bZ$
and $x$ in the domain of $\t_{n i}.$ Thus, for an element $x \in \ar$ 
that is analytic of exponential type with respect to $\t$ 
in the sense of \cite[Definition 4.1]{Ha}, it follows from complex function 
theory (see e.\ g.\ \cite[Lemma 5.5]{Wo1}) that $\eta(\t_{z}(x) )= \eta(x)$ 
for all $z \in \bC.$ Now, the set of such elements in $\ar$ is easily seen 
to be invariant under each $\t_{t}, t \in \bR,$ and dense in $\ar$ 
(by the proof of \cite[Lemma 4.2]{Ha}).
Hence  Corollary 1.22 
of \cite{Ku2} says that this set is a core for any $\t_{z}, z \in \bC.$ 
Thus, for any $z \in \bC,$ we have $\eta(\t_{z}(x) )= \eta(x)$ for all 
$x$ in the domain of $\t_{z}.$ In particular, choosing $z = i/2,$ we get 
$\eta \tau_{i/2} = \eta$ on $\pi(A)$ as asserted.)  
Using \cite[Theorem 5.6]{KD}, we then get
\[ \eta R \pi(a) = \eta \tau_{i/2} \pi(S(a)) = \eta \pi(S(a)), \]
for all $a \in A$. This gives
\[(\eta \otimes \eta R ) \dar (\pi (a)) = (\eta \pi \otimes \eta R \pi )
(\D(a)) =  (\eta \pi \otimes \eta  \pi S )  (\D(a))\]
\[= \eta \pi (m (\iota \otimes S) \D (a)) = \e(a), \]
where ${m:A\otimes A\to A}$ is the linearization of the
multiplication of~$A$. Thus, $(\eta \otimes \eta R ) \dar( \pi
(a)) = \e(a)$, for all $a \in A$, as required, and
Condition~(1) holds.

Now suppose Condition (5) holds and let $\eta$ be a unit for
$\ar^{*}$. For $a \in A$ and $\r \in \ar^{*}$ we then have
\[ \r(\pi(a)) = (\eta \r)(\pi(a)) = \r ( (\eta \otimes \iota)
\dar(\pi(a)) ). \] Since $\ar^{*}$ separates $\ar$
 we get $(\eta \otimes \iota) \dar(\pi(a))
= \pi (a) $, for all $a \in A$. In the same way we also get
$(\iota \otimes \eta) \dar(\pi(a)) = \pi (a) $, for all $a \in
A$. From the uniqueness property of the co-unit we can then
conclude that $\eta \pi = \e$.
 Since $\eta$ is bounded by assumption, it follows that $\e$ is
 bounded (with
 respect to the norm on $A$ inherited from the one on $\pi(A)$). Hence
 Condition~(1) holds.

Suppose Condition (2) holds, so that there exists a net
$(v_i)$ of unit vectors in $H$ such that $\lim_i\|
W(v_i\otimes v)-v_i\otimes v\|_2 =0$, for all $v\in H$. Define
$\omega_{i} \in M_*$ to be the restriction of $ \omega_{v_i}$
to $M$.
 Then, for all  $v \in H$ and
all $x$ in the unit ball of $M$, we have
 \[  | (\omega_{i} \omega_{v})(x) - \omega_{v}(x) | =  | (\omega_{i}
 \otimes \omega_{v}) ( W^{*} (1 \otimes x ) W) - \omega_{v}(x) |   \]
 \[ =| ( W^{*} (1 \otimes x ) W (v_{i}\otimes
 v) - (1 \otimes x ) (v_{i} \otimes v) \, , \, v_{i} \otimes v )   |  \]
 \[ =  |( (1 \otimes x ) W (v_{i}\otimes
 v) \, , \,  W (v_{i}\otimes v) )
 - ( (1 \otimes x ) (v_{i} \otimes v) \, , \, v_{i} \otimes v ) | \]
 \[   =  | ( (1 \otimes x ) (W (v_{i}\otimes
 v) - v_{i}\otimes v)  ,   W (v_{i}\otimes v) )
 + ( (1 \otimes x ) (v_{i} \otimes v) ,  W (v_{i}\otimes
 v) - v_{i}\otimes v)   | \]
 \[ \leq 2 \| v \| \, \| W (v_{i}\otimes
 v) - v_{i}\otimes v \|_{2}. \]

Hence, \[ \| \omega_{i} \omega_{v} - \omega_{v} \|\leq 2 \| v
\| \, \| W (v_{i}\otimes
 v) - v_{i}\otimes v \|_{2}  \rightarrow 0. \]
 Since $M$ is in standard form, every normal state on $M$
 is equal to (the restriction of) $\om_v$, for some unit vector $v\in H.$
 It follows therefore from our calculations that
 $(\omega_{i})$ is a bounded left
 approximate unit for $M_{*}$. Hence, Condition~(2) implies Condition~(6).

 To see that (6) $\implies$ (7) and  (7) $\implies$ (8), let us just
 remark that if $(\om_i)$ is a bounded left approximate unit
 for $M_*$, and we set
 $\omega_{i}^{o} = \omega_{i} R  \in M_{*}$, then,
 using the fact that $\chi (R \bar{\otimes} R) \dar = \dar R$, it is
 straightforward to
 check that   $(\omega_{i}^{o})$ is a bounded right
 approximate unit for $M_{*}$. The map $\chi$ is, of course, the
 flip map on ${M\bar\otimes M}$. It is then easily seen that
 $(\omega_{i} +
 \omega_{i}^{o} - \omega_{i}^{o}\omega_{i})$ is a
 bounded two-sided  approximate unit for $M_{*}$.

Finally assume Condition~(8) holds and that $(\om_{i})$ is
bounded two-sided approximate unit for $M_{*}$. By going to a
subnet of $(\om_i)$, if necessary, we may suppose that
$(\om_i)$ converges in the weak$^*$-topology in $M^*$ to an element $\om$.
We use the same symbol to denote an element in $M^{*}$ and its 
restriction to $\ar$.
Let  $x \in \pi(A)$. Then, for all $\om' \in M_{*}$, we have 
\[ \om' ( (\om \otimes \iota ) \dar (x) ) = \om ( (\iota 
\otimes \om') \dar (x) ) = \lim_{i} \ \om_{i} ((\iota 
\otimes \om') \dar (x) ) \]
\[ = \lim_{i} \ \om_{i} ((\iota 
\bar\otimes \om') \dar (x) ) = \lim_{i} \ \om_{i} \om'  (x) = \om' (x).\]
Since the set $M_*$ separates the elements
of $M$, it follows that ${(\om \otimes \iota) \dar(x)} = x$. 
Similarly, we get $ {(\iota \otimes \om) \dar(x)} = x$.
 Hence, for all $a\in A$, $(\om\otimes \iota) \dar(\pi(a)) = (\iota
 \otimes \om) \dar(\pi(a))=\pi(a)$. From the uniqueness
 property of the co-unit, we conclude that $\om \pi = \e$.
 Since $\om$ is norm-bounded, it follows that $\e$ is
 norm-bounded also. Hence, Condition~(8) implies~(1). This
 completes the proof of the theorem.~\qed

To each algebraic quantum group $\AA$ one may construct
a unique universal C$^{*}$-algebraic quantum group $(A_{u},\Delta_{u})$ (see
\cite{Ku1}). Co-amenability of $\AA$ may be seen to be equivalent to the
fact
that the canonical homomorphism from $A_{u}$ onto $\ar$ is
injective (see \cite{BMT} for the compact case). We will return to this aspect
of co-amenability in a subsequent work \cite{BMT2}.
 
In the case that the \aqg\ $\AA$ is of compact type, we can
prove some results that make explicit use of the existence of
the unit in~$A$. In this case we can choose a unique
left-invariant unital linear functional $\f$ on~$A$. This is
the left Haar integral of $A$ and it is also right invariant. We
refer to $\f$ as the {\em Haar state} of~$A$. If $(H,\Lambda)$
is a GNS~pair associated to $\f$, then the restriction
$\f_{\mathrm{r}}$ of the state  $\om_{\Lambda(1)}$ to $\ar$ is
a left and right invariant state for the co-multiplication
${\dar:\ar\to \ar\otimes \ar}$.

The following result generalizes Lemma~10.2 of~\cite{MT}.

\begin{lem}
Let $\AA$ be an \aqg\ of compact type and let $B$ be a
C*-algebra admitting a faithful state (for example, let $B$ be
separable). Let $\ta\colon A\to B$ be a unital $*$-linear map
that is either multiplicative or anti-multiplicative. Then the
linear map, ${\ta'\colon A\to A\otimes B}$, ${a\mapsto
(\iota\otimes \ta)\D(a)}$, is isometric.

\end{lem}

\demo We shall prove the result in the multiplicative case
only---the proof in the anti-multiplicative case is similar.
We identify $A$ as a $*$-subalgebra of $\ar $. Let $\t$ be a
faithful state on~$B$. Since the Haar state $\f_{\mathrm{r}}$
on $\ar $ is faithful, the state ${\f_{\mathrm{r}}\otimes \t}$
on~${\ar \otimes B}$ is faithful. Hence, by two applications
of \cite[Theorem~10.1]{MT} (to $\f_{\mathrm{r}}\otimes \t$ and
then to $\f_{\mathrm{r}}$), if $a\in A$, we get
$\n{\ta'(a)}^2=\n{\ta'(a)^*\ta'(a)}=\lim[(\f_{\mathrm{r}}\otimes
\t)(\ta'(a^*a)^n)]^{1/n}=\lim[(\f_{\mathrm{r}}\otimes
\t\ta)\D((a^*a)^n)]^{1/n}=\lim[\t\ta(1)\f_{\mathrm{r}}((a^*a)^n)]^{1/n}=
\n{a^*a}=\n a^2$. Thus, $\ta'$ is isometric, as required.~\qed

\begin{thm} \label{thm: finite dimensional reps}
Let $\AA$ be an \aqg\ of compact type and suppose there exists
a non-zero continuous $*$-homomorphism $\ta$ from $A$ onto a
finite-dimensional C*-algebra~$B$. Suppose also that for the
antipode $S$ of $A$ we have  $\ta S^2=\ta$. Then $\AA$ is
co-amenable.
\end{thm}

\demo Let $1_B$ denote the unit of $B$ and let the map
${\e_B\colon A\to B}$ be defined by setting
$\e_B(a)=\e(a)1_B$. Let ${m_B\colon B\otimes B\to B}$ and
${m\colon A\otimes A\to A}$ be the multiplications. Since
$\ta$ is multiplicative, $m_B\circ(\ta\otimes \ta)=\ta\circ
m$. Using this and the fact that $m(\iota\otimes
S)\D(a)=\e(a)1$, for all $a\in A$, we get
$$\e_B=m_B\circ(\ta\otimes \iota)\circ[(\iota\otimes \ta
S)\circ\D].$$ This is a product of three continuous linear
maps. The map $m_B$ is continuous, since $B$ is
finite-dimensional; ${\ta\otimes \iota}$ is continuous, since
$\ta$ is; finally, the map $(\iota \otimes \ta S)\circ\D$ is
continuous by the preceding lemma, since $\ta S$ is obviously
unital, $*$-linear and  anti-multiplicative. Hence, $\e_B$ is
continuous. It follows immediately that $\e$ is continuous and
therefore that $\AA$ is co-amenable.~\qed

The assertion (4) implies (1) in Theorem~\ref{thm: net
equivalence} may be rephrased as saying that  $\AA$ is
co-amenable whenever there exists a non-zero continuous
complex valued homomorphism on $A$. Note that we etablished 
in our proof of this fact that
such a homomorphism is always $S^2$-invariant. On the other
hand, we don't know whether the $S^{2}$-invariance assumption
in Theorem~\ref{thm: finite dimensional reps} is redundant.

\begin{cor} \label{thm: S^2=1}
Let $\AA$ be an \aqg\ of compact type and suppose that its
analytic extension~$\ar $ is of Type~I, as a C*-algebra.
Suppose also that for the antipode $S$ of $A$ we have
$S^2=\iota_A$. Then $\AA$ is co-amenable.
\end{cor}

\demo Since $\ar $ is unital, it admits a maximal ideal $I$.
The quotient algebra $B=\ar /I$ is a C*-algebra of Type~I and
is both unital and simple. Therefore, it  is finite
dimensional. Hence, the restriction of the quotient map is a
non-zero continuous $*$-homomorphism $\ta$ from $A$ onto a
finite-dimensional C*-algebra, namely~$B$. From the assumption
that $S^2=\iota_A$, the existence of $\ta$ implies
co-amenability of $\AA$, by the theorem above.~\qed

It would be interesting to know whether Corollary~\ref{thm:
S^2=1} remains true in the case that $S^2\ne \iota_A$; that
is, where the Haar state of $\AA$ is not tracial
(\cite{BS,Wo1}).

\vspace{2ex} We now define a notion ``dual'' to
co-amenability, namely amenability. As pointed out in the
introduction, this notion is due to D. Voiculescu in the Kac
algebra case~\cite{Vo} (see also \cite{ES,Ruan}).

Let $\AA$ be an \aqg\ with von~Neumann algebra~$M$. A {\em
right-invariant mean} for $(A, \D)$ is a state $m$ on $M$ such
that
\[m((\iota\bar{\otimes}\omega) \dar(x))=\omega (1)m(x),\]
for all $x\in M$ and $\omega\in M_*$. A {\em left-invariant
mean} is defined analogously, but we shall have no need for
this concept in this paper. We say that $\AA$ is {\em
amenable} if $(A, \D)$ admits a right-invariant mean. Using
the existence of the anti-unitary antipode $R$ on $(M, \dar)$
(\cite{KD,KV}), this is easily seen to be equivalent to
requiring that $(A, \D)$ admits a left-invariant mean.

\begin{example} \label{group example continued} \rm The
relation to amenability in the classical case of a discrete
group $\Gamma$ is worth considering in some little detail.
Recall from Example~\ref{group example} that the group algebra
$\AA=(\bC(\Gamma),\D)$ is an algebraic quantum group with dual
$(\hat A,\hat\D)=(\fs(\Gamma),\hat\D)$ and that the map,
${\hat\pi:\hat A\to B(\ell^2(\Gamma))}$, extends to a
$*$-isomorphism from $\ell^\infty(\Gamma)$ onto $\hat M$ and
that this $*$-isomorphism is just the usual representation by
multiplication operators. For $x\in \Gamma$, let
${R_x:\ell^\infty(\Gamma)\to \ell^\infty(\Gamma)}$ be the
right translation operator given by $R_xf(y)=f(y x)$, for all
$y\in \Gamma$. We claim that, for all $f\in
\ell^\infty(\Gamma)$,
\begin{equation}\hat\pi R_x(f)=
(\iota\bar\otimes\om_{e_x})\hat\dar\hat\pi(f),\label{weak
equation}
\end{equation}
where $e_x$ is defined as in Example~\ref{group example}.
Since both sides of this equation belong to the algebra $\hat
M=\hat\pi(\ell^\infty(\Gamma))$ and this is diagonal with
respect to the orthonormal basis $(e_x)_{x\in \Gamma}$, to see
the equality holds we need only show, for all $y\in \Gamma$,
\[\om_{e_y}(\hat\pi R_x(f))=
\om_{e_y}((\iota\bar\otimes\om_{e_x})\hat\dar\hat\pi(f)).\]
Now the functional $\om_{e_y}\hat\pi$ on $\ell^\infty(\Gamma)$
is easily verified to be just the operation of evaluation at
$y$, so we need only show \begin{equation} \label{second weak
equation} f(y x)=\ip{W(\hat\pi(f)\otimes 1) W^*(e_y\otimes
e_x)}{(e_y\otimes e_x)}.\end{equation} However, direct
computation shows that ${W^*(e_y\otimes e_x)}={e_{yx}\otimes
e_x}$. Hence,
\[\ip{W(\hat\pi(f)\otimes 1)W^*(e_y\otimes e_x)}{e_y\otimes
e_x}=\ip{(f\otimes 1)(e_{yx}\otimes e_x)}{e_{yx}\otimes
e_x}\]\[=\ip{fe_{yx}}{e_{yx}}=f(y x).\] This shows that
Equation~(\ref{second weak equation}) holds and therefore
Equation~(\ref{weak equation}) also holds.

Suppose now $m$ is a state on $\hat M$ and let $\tilde m$ be
the state on $\ell^\infty(\Gamma)$ determined by $\tilde
m\circ \hat\pi=m$. Using the fact that $\hat M$ is diagonal
with respect to the orthonormal basis $(e_x)_{x\in \Gamma}$,
it is easily checked that $m$ is a right-invariant mean on
$\hat M$ if, and only if,
$m((\iota\bar\otimes\om_{e_x})\hat\dar\hat\pi(f))=m\hat\pi
(f)$, for all $x\in \Gamma$ and $f\in \ell^\infty(\Gamma)$.
Using Equation~(\ref{weak equation}) this translates into the
condition $m\hat\pi(R_xf)=m\hat\pi(f)$; that is, $\tilde
m(R_xf)=\tilde m(f)$. Thus, $(\fs(\Gamma),\hat\D)=(\hat A,\hat
\D)$ is amenable if, and only if, the group $\Gamma$ is
amenable.

Note that $(\bC(\Gamma), \D)=\AA$ is always amenable, since
$\om_{\Lambda(1)}$ provides a right-invariant mean on~$M$, as
is easily verified.~\qed \end{example}

Amenability is automatic for an \aqg\ $\AA$ of compact type.
If $\f$ is Haar state of $A$ and $(H,\Lambda)$ is a GNS~pair
for $\f$, then the vector state $m=\om_\Lambda(1)$ is easily
seen to define a right-invariant mean on~$M$.

\begin{thm}
If the algebraic quantum group $(A, \D )$  is co-amenable,
then the dual algebraic quantum group $(\hat{A}, \hat{\D})$ is
amenable.
\end{thm}

\demo Let $(H,\LL)$ be a GNS~pair, and $W$ the
corresponding multiplicative unitary, associated to $(A, \D)$. 
Assume that $(A, \D )$  is co-amenable and let $\e_{\mathrm{r}}$ be a state on $\ar$ 
such that $(\e_{\mathrm{r}}\otimes\iota)(W)=1,$ according to Theorem \ref{thm: net 
equivalence}. Let $\e'$ denote a state extension of 
$\e_{\mathrm{r}}$ to $B(H).$ 
Then the restriction $m$ of $\e'$ to $\hat M$ is a
right-invariant mean for $(\hat{A}, \hat{\D}).$

Indeed let $v$ be a unit vector in~$H$. To see that the
restriction is right-invariant, we need only show that
\[{m(\iota\bar\otimes
\om_v)(\hat\dar(x)}={m(\iota\bar\otimes
\om_v)(W(x\otimes 1)W^*)}=m(x), \ x \in \hat M . \] Observe that
$\e'(\iota\bar\otimes \om_v)$ is a state on 
${B(H\otimes H)}=B(H)\bar\otimes B(H)$ and
$${\e'(\iota\bar\otimes \om_v)(W)} = \e_{\mathrm{r}}  (\iota \otimes 
\om_{v})(W) = {\om_v(\e_{\mathrm{r}}\otimes
\iota)(W)}=\om_v(1)=1.$$ Hence $\e'(\iota\bar\otimes \om_v)$ is 
multiplicative at $W$ and $W^*$ (see our remark at the end of section 2). It follows that
$${m(\iota\bar{\otimes} \om_v)(W(x\otimes
1)W^*)}= {\e'(\iota\bar{\otimes} \om_v)(W(x\otimes
1)W^*)}$$
$$={\e'(\iota\bar\otimes \om_v)(x\otimes 1)}=
{\e'(x)\om_v(1)}=m(x)$$ for all $x \in \hat M,$  as required.~\qed

The preceding theorem raises the question as to whether its
converse holds; that is, if the algebraic quantum group $(A,
\D )$  is such that its dual $(\hat{A}, \hat{\D})$ is
amenable, is $\AA$ co-amenable? Recall from Example~\ref{group
example} that if $\AA$ is the \aqg\ $(\bC(\Gamma),\D)$, where
$\Gamma$ is a discrete group, its dual is $(\hat
A,\hat\D)=(\fs(\Gamma),\hat\D)$. Hence, $(\hat A,\hat \D)$ is
amenable if, and only if, $\Gamma$ is amenable, as we saw in
Example~\ref{group example continued}. On the other hand,
 as mentioned before, it is well known that $\AA$
is co-amenable if, and only if, $\Gamma$ is amenable. Thus, in
this case, the converse of the preceding theorem holds.  In
the more general case that $\AA$ is of compact type and the
Haar functional is tracial, the above question may also be
answered positively as may be deduced from Ruan's
main result (Theorem 4.5) in \cite{Ruan}. We will discuss this fact
and related matters in a subsequent paper (\cite{BMT2}).
 
\vspace{2ex}
Recall that amenability of C$^{*}$-algebras and of von Neumann algebras
may be described
by several equivalent formulations (see \cite{P} for an
overview of these and for references to the literature). The
most commonly used terminology is nuclearity for C*-algebras
and injectivity for von~Neumann algebras. The following result
may be seen as a quantum group counterpart of the well known
result that the group von~Neumann algebra of a (locally
compact) group is injective whenever the group is amenable.

We shall need some easy preliminaries on
generalised limits for our next theorem. Let $(I,\le)$ be a
directed pair; that is, $I$ is a non-empty set and $\le$ is a
reflexive, transitive relation on~$I$ that is directed upwards
in the sense that for all $i,j\in I$, there exists $k\in I$
such that $i,j\le k$. Then there is a norm-decreasing linear
functional $\om$ on $\ell^\infty(I)$ such that $\om(x)=\lim_i
x_i$, for each convergent sequence $x$ in $\ell^\infty(I)$. We
call any such functional $\om$ a {\it generalised limit
functional} for $(I,\le)$. Obviously, $\om$ is a state of
$\ell^\infty(I)$. It is elementary to see that a generalised
limit functional $\om$ exists. First, define $\om$ as a
norm-decreasing linear functional in the obvious way on the
linear subspace of $\ell^\infty(I)$ consisting of elements $x$
that are eventually constant in the sense that there is a
scalar $\l$ and an element $i\in I$ for which $x(j)=\l$, for
all $j\ge i$. Then use the Hahn--Banach theorem to extend
$\om$ to a norm-decreasing linear functional
on~$\ell^\infty(I)$. (This is  a weakening of the usual
concept of a Banach limit on the positive integers.)

\begin{thm} \label{thm:injectivity}
Let $\AA$ be a co-amenable \aqg. Then its von~Neumann algebra
$M$ is  injective.
\end{thm}

\demo Of course, we have to show that there is an norm-bounded
idempotent operator $E$ on~$B(H)$ with range $M$. As usual,
$(H,\LL)$ is a GNS~pair associated to a left Haar integral
$\f$ on $A$ and $W$ is the corresponding multiplicative
unitary. By the proof of Theorem~\ref{thm: net
equivalence}~(7) and co-amenability of $\AA$, there is a net
$(\om_{i})_{i\in I}$ of normal states on $M$ which is a right
approximate unit for $M_{*}$. As $M$ is in the standard
representation, we may, and do, choose a net $(v_i)_{i\in I}$
of unit vectors in~$H$ such that, $\om_{i}$ is the restriction
of $\om_{v_i}$ to $M$. Now choose a generalised limit
functional on $\ell^\infty(I)$ and denote its value at a
bounded net $(x_i)_i$ in $\bC$ by $\Lim_i x_i$.

Let $x\in B(H)$. If $v,w\in H$, then
${((\om_{v,w}\bar{\otimes}\om_{v_i})(W^*(1\otimes x)W))_i}$ is
a bounded net. Set $\eta_x(v,w)=\Lim_i
{(\om_{v,w}\bar{\otimes}\om_{v_i})(W^*(1\otimes x)W)}$. It is
easily verified that this defines a sesquilinear form and that
$\m{\eta_x(v,w)}\le \n x\n v \n w$. Hence, there is a unique
linear operator $E(x)$ such that $\ip{E(x)v}w=\eta_x(v,w)$, for
all $v,w\in H$. It follows that $\n{E(x)}\le \n x$. The map,
${E\colon B(H)\to B(H)}$, ${x\mapsto E(x)}$, is obviously
linear.

We shall show next that $E(B(H))\subseteq M$. 
Suppose that $x\in B(H)$ and $E(x)\notin M$. Then there
exists a unitary $U\in M'$ such that $U^*E(x)U\ne E(x)$.
Hence, there exists an element $v\in H$ such that
$\ip{U^*E(x)U v}v\ne \ip{E(x)v}v$. That is,
$\om_{U(v)}(E(x))\ne \om_v(E(x))$. Set $\t=\om_{U(v)}-\om_v$.
Then $\t(E(x))\ne 0$. However, we clearly have  $\t(M)=0$.
Hence, for all $i\in I$,
${(\t\bar{\otimes}\om_{v_i})(W^*(1\otimes x)W)=0}$, since
$W^*(1\otimes x)W\in {M\bar{\otimes} B(H)}$. It follows that
${(\om_{U(v)}\bar{\otimes}\om_{v_i})(W^*(1\otimes x)W)}=
{(\om_v\bar{\otimes}\om_{v_i})(W^*(1\otimes x)W)}$. From this
we get  \[{\Lim_i
(\om_{U(v)}\bar{\otimes}\om_{v_i})(W^*(1\otimes x)W)}= {\Lim_i
(\om_v\bar{\otimes}\om_{v_i})(W^*(1\otimes x)W)};\] that is,
$\ip{E(x)U(v)}{U(v)}=\ip{E(x)v}v$ and therefore $\t(E(x))=0$.
This is a contradiction and to avoid it we must have $E(x)\in
M$.

To complete the proof we need only show now that $E(x)=x$, for
all ${x\in M}$. We have, for each element $v$ in $H$,
\[\ip{E(x)v}v={\Lim_i (\om_v\bar{\otimes}\om_{v_i})(W^*(1\otimes
x)W)} = \Lim_i(\om_v \om_{i}) (x)=\om_v(x).\] Consequently,
$\ip{E(x)v}v=\ip{xv}v$, for all $v\in H$, and therefore
$E(x)=x$, as required.~\qed

After some work, one may see that the above result can  be
deduced from \cite[Proposition 3.10]{Ng1}, which itself may be seen
as a generalization of \cite[Proposition 5.6 ]{Bl}. Both of
these results deals with nuclearity (of crossed products in
the context of Hopf C*-algebras and in the context of  regular
multiplicative unitaries, respectively). Injectivity of
crossed products in the Kac algebra case has been considered
in \cite[Section~3]{ES}. On the other hand, the converse to
Theorem~\ref{thm:injectivity} is known to hold in the compact tracial
case, as may be obtained from \cite[Theorem 4.5]{Ruan} (see 
\cite{BMT2} for a simplified proof of this fact). To illustrate the 
concepts we include here  a direct proof of a related result on ``dual'' 
amenability.

\begin{thm} \label{thm: injectivity implies dual amenability}
Suppose that $\AA$ is an algebraic quantum group $(A,\D )$ of
compact type and that its Haar state~$\f$ is tracial. If
$M$ is injective, then the dual \aqg\ $(\hat{A},\hat{\D})$ is
amenable.
\end{thm}

\demo It is clear that $\om_w$, where $w=\Lambda(1)$, is a
normal state on $M$ and that $\om_w\circ\pi=\f$, so that
$\om_w$ is the ``analytic extension'' of $\f$ to $M$. The
hypothesis gives the existence of a unital norm-decreasing
positive idempotent linear operator ${E:B(H)\to B(H)}$ such
that $E(B(H))=M$. Define a norm-decreasing linear functional
$m$ on $\hat M$ by setting $m(x)=\om_w(E(x))$, for all $x\in
\hat M$. Since $\om_w$ and $E$ are unital, so is $m$ and
therefore $1=m(1)=\n m$. Hence, $m$ is a state. We are going
to prove now that $m$ is a right-invariant mean on $\hat M$.

To show this, we need only show, for all unit vectors $v\in H$
and all $x\in \hat M$,
\begin{equation}m((\iota\bar\otimes\om_v)(\hat\dar(x)))
=m(x). \label{m-equation}\end{equation}  Using the fact that
$\Lambda(A)$ is dense in $H$, we may suppose that
$v=\Lambda(a)$; hence, $\f(a^*a)=1$. Since  $\hat{A}_r$ acts 
non-degenerately
on $H$ and $\hat \pi (\hat A )$ is dense in $\hat{A}_r,$ there
is a norm-bounded net $(e_i)$ in $\hat A$ such that $((\hat
{\pi}(e_i)(v))$ converges to $v$. Hence, the net of vector
states $(\omega_{\hat{\pi}(e_i )v})$ on~$\hat M$ converges
to $\omega_v$ in norm, and therefore
\[\|(\iota\bar{\otimes}\omega_v )\hat \D_{\mathrm{r}}
(x)-(\iota\bar{\otimes}\omega_{\hat{\pi}(e_i )v})\hat
\D_{\mathrm{r}} (x)\| \leq\|\omega_v -\omega_{\hat{\pi}(e_i
)v}\|\|x\|\rightarrow 0.\] Recall that $W\in M(\ar \otimes
\hat \ar )$, from which it follows that $Z_i=(\pi
(1)\otimes\hat{\pi}(e_i^* ))W$ and $Z_i^*=W^*(\pi
(1)\otimes\hat{\pi}(e_i ))$ both belong to
$\ar\otimes\hat{\ar}$, for all indices~$i$. Now write $\D a
=\sum_j a_j\otimes b_j$, for some elements $a_j,b_j\in A$.
Then
\[m((\iota\bar{\otimes}\omega_v )\hat \D_{\mathrm{r}} (x))=
\lim_i \om_w\circ E ((\iota\bar{\otimes}\omega_{\hat{\pi}(e_i
)v} )\hat \D_{\mathrm{r}} (x))\]
\[=\lim_i \om_w\circ E\circ(\iota\bar{\otimes}\omega_v)
((\pi (1) \otimes\hat{\pi}(e_i^* ))W(x\otimes 1)W^* (\pi (1)
\otimes\hat{\pi}(e_i )))\]
\[=\lim_i (\om_w E\otimes \om_v)(Z_i(x\otimes 1)Z_i^*)=
\lim_i (\om_w \otimes \om_v)(Z_i(E(x)\otimes 1)Z^*_i).\] For
the last equality we are using the fact that $E(c y d)=c E(y)
d$, for all $c,d\in M$ and $y\in B(H)$. Hence,
\[m((\iota\bar{\otimes}\omega_v )\hat
\D_{\mathrm{r}} (x))=\lim_i ((W(E(x)\otimes 1)W^*
(\LL(1)\otimes \hat{\pi} (e_i )v),\LL(1)\otimes \hat{\pi} (e_i
)v)\]
\[=(W(E(x)\otimes 1)W^* (\LL(1)\otimes \Lambda (a)),
\LL(1)\otimes \Lambda (a))\]
\[=((E(x)\otimes 1)W^* (\LL(1)\otimes \Lambda (a)),
W^* (\LL(1)\otimes \Lambda (a)))\]
\[=((E(x)\otimes 1)(\LL\otimes \LL)(\D(a)),
(\LL\otimes \LL)(\D(a))\]
\[=\sum_{jk} ((E(x)\otimes 1)(\pi(a_j)\LL(1)\otimes \pi(b_j)\LL(1)),
\pi(a_k)\LL(1)\otimes \pi(b_k)\LL(1))\]
\[=\sum_{jk}\om_w (\pi(a_k)^* E(x)\pi(a_j) )\varphi (b_k^*
b_j ).\]

Thus, to show $m$ is right-invariant, it is sufficient to show
that
\[\sum_{jk}\om_w (\pi(a_k)^* E(x)\pi(a_j) )\varphi (b_k^*
b_j )=m(x)=\om_w(E(x)),\] for all $x\in \hat M$. Using the
fact that $\om_w$ is weakly continuous, it therefore suffices
to show that \[\sum_{jk}\om_w (\pi(a_k)^* \pi(b)\pi(a_j)
)\varphi (b_k^* b_j )=\om_w(\pi(b)),\] for all $b\in A$; that
is,
\[\sum_{jk}\f(a_k^* ba_j )\varphi (b_k^* b_j )=\f(b).\]
However, the left side of this equation is
$(\f\otimes\f)(\D(a)^*(b\otimes 1)\D(a))$ and this is equal to
the right side of the equation, since $\f(a^*a)=1$ and $\f$ is
tracial. Consequently, we have shown
$m((\iota\bar{\otimes}\omega_v )\hat \D_{\mathrm{r}}
(x))=m(x)$, for all $x\in \hat M$. Therefore,
$(\hat{A},\hat{\D})$ is amenable.~\qed

Suppose that $\AA$ is an \aqg\ of compact type and its Haar
state is tracial. Since nuclearity of $\ar$ implies
injectivity of $M$, it follows from Theorem~\ref{thm:
injectivity implies dual amenability} that nuclearity of $\ar$
implies amenability of $(\hat A,\hat \D)$.

Note also that the tracial assumption in Theorem~\ref{thm: injectivity
implies dual
amenability} is equivalent to assuming that $S^2=\iota$ (\cite{BS, Wo1}).

\section{Co-amenability and modular properties}

In this section we investigate the modular properties of a
co-amenable algebraic quantum group $(A, \D )$ of compact
type. The unital Haar functional $\varphi$ of $(A, \D )$ is a
KMS-state when extended to~$\ar $. In the case that $\AA$ is
co-amenable, we show that the modular group can be given a
description in terms of the multiplicative unitary of~$\AA$.

We shall make use of the existence of a certain family
$(f_z)_{z\in \bC}$ of functionals on~$A$; these functionals
are quite particular to the compact quantum group case. Observe here
that the C$^{*}$-algebraic compact quantum group  $(\ar,\dar)$ has
clearly $\AA$ as its canonical dense
Hopf $*$-algebra (see \cite[Appendix]{BMT}), and that we may therefore
use  Woronowicz's results from \cite{Wo1,Wo}. Before
stating the properties of these functionals, let us recall
that an entire function~${g\colon\bC\to \bC}$ has {\em
exponential growth on the right half plane} if there exist
real numbers $M$ and $r$, with $M>0$, such that $\m{g(z)}\le
Me^{r\Re(z)}$, for all $z\in \bC$ for which $\Re(z)\ge 0$.

There exists a unique family $(f_z)_{z\in{\bf C}}$ of unital
multiplicative linear functionals on~$A$, that we shall call
the {\em modular functionals}, satisfying the following
conditions:

\begin{enumerate}
\item For each element $a\in A$, the map $z\mapsto f_z (a)$ is a
entire function of exponential growth on the right half-plane;
\item $f_0 =\varepsilon$ and $f_{z+w}  =(f_z\otimes f_w )\D$,
for all $z,w\in \bC$;
\item $f_z (S(a))=f_{-z}(a)$ and $\overline{f_z (a^*)}=
f_{-\overline{z}}(a)$, for all $a\in A$ and $z\in \bC$;
\item $S^2 (a)=f_{-1}* a* f_1$, for all $a\in A$;
\item The unique homomorphism ${\r\colon A\to A}$ for which
$\varphi (a b)=\varphi (b \r(a))$, for all $a,b\in A$, is
given by $\r(a)=f_1*a*f_1$.
\end{enumerate}

Here $\t*a*\t' =(\t'\otimes\iota\otimes\t)(\D\otimes\iota)\D
(a)$, for all $a\in A$ and all linear functionals $\t, \t'$ on
$A$. Note in particular that it follows from Condition~3 that
$f_{it}$ is a unital $*$-homomorphism from $A$ to $\bf C$,
for all real numbers~$t$. Hence, for each element $t\in \bR$,
one can define a unital $*$-automorphism $\sigma_t$ of
$\pi(A)$ by setting
\[\sigma_t (\pi (a))=\pi (f_{-it}*a*f_{-it}),\]
for all $a\in A$. Since $\varphi_{\mathrm{r}}\sigma_t
=\varphi_{\mathrm{r}}$ and $\| a\|^2 =\lim_n
\varphi_{\mathrm{r}} ((a^* a)^n )^{1/n}$, for all ${a\in \ar
}$, the $*$-homomorphisms $(\sigma_t)_{t\in {\bf R}}$ are
isometric and can therefore be extended to a one-parameter
group of $*$-automorphisms of $\ar $, also denoted by
$(\sigma_t)_{t\in{\bf R}}$. Moreover, we have a unitary
implementation. This is easily seen by noting first that the
linear map from $\Lambda(A)$ to itself defined by sending
$\Lambda(a)$ onto $\Lambda(f_{-it}*a*f_{-it})$ is an isometric
bijection. Its extension to $H$ is a unitary $V_t$. One easily
verifies then that $\sigma_t (x)=V_t x V_t ^*$, for all $x\in
\pi(A)$ and therefore for all $x\in \ar $.

However, if we now suppose that $\AA$ is co-amenable, more can
be said about the unitaries $\{V_t\}_{t\in{\bf R}}$. In this
case it follows from \cite[Corollary 3.7]{BMT} that all positive
functionals on~$A$ are bounded  with respect to the norm on
$A$ (the one induced from the norm on~$\ar $). Hence, there is
an extension $\tilde{f}_{it}$ of $f_{it}$ to~$\ar $; more
precisely, there is a unique $*$-homomorphism $\tilde{f}_{it}$
on $\ar $ such that $\tilde{f}_{it} \pi (a)=f_{it}(a)$, for
all $a\in A$. The map ${\tilde{f}_{it}\otimes\iota :\ar
\otimes B_0 (H)\rightarrow B_0 (H)}$ is clearly a
non-degenerate $^*$-homomorphism and therefore has a unique
extension to a $*$-homomorphism ${\tilde{f}_{it}\otimes\iota
:M(\ar \otimes B_0 (H) )\rightarrow M(B_0 (H))=B(H)}$. Since
$W\in M(\ar \otimes B_0(H) )$,  we may define an element
$F_t\in B(H)$ by setting $F_t
={(\tilde{f}_{it}\otimes\iota)(W)}$. Now observe that
$\tilde{f}_{i(s+t)}=(\tilde{f}_{is}\otimes\tilde{f}_{it})\dar
$, for all $s,t\in\bf R$. This result follows by extending the
formula $f_{is +it }={(f_{is}\otimes f_{it})\Delta}$ from $A$
to~$\ar $. Using the equation ${(\dar \otimes\iota
)W}={W_{13}W_{23}}$ we get
\[F_{s+t}=(\tilde{f}_{i(s+t)}\otimes\iota)(W)=
((\tilde{f}_{is}\otimes\tilde{f}_{it})\dar \otimes\iota )(W)
=(\tilde{f}_{is}\otimes\tilde{f}_{it}\otimes\iota )(\dar
\otimes\iota )(W)\]
\[=(\tilde{f}_{is}\otimes\tilde{f}_{it}\otimes\iota )(W_{13}W_{23} )
=(\tilde{f}_{is}\otimes\iota )(W)(\tilde{f}_{it}\otimes\iota
)(W) =F_s F_t.\] Since we have seen that $F_0
=(\tilde{\varepsilon}\otimes\iota )(W) =1$, where the
functional $\tilde{\varepsilon} :\ar \rightarrow\bf C$ is the
extension of the co-unit of $(A, \D )$, it is  clear that the
operators $F_t$ are invertible in $B(H)$. In fact, $F_t$ is a
unitary, since
\[F_t ^* F_t =((\tilde{f}_{it}\otimes\iota)(W))^*
(\tilde{f}_{it}\otimes\iota)(W)
=(\tilde{f}_{it}\otimes\iota)(W^*
)(\tilde{f}_{it}\otimes\iota)(W)\]
\[=(\tilde{f}_{it}\otimes\iota)(W^* W )=
(\tilde{f}_{it}\otimes\iota)(1)=1.\] Hence, the map,
${t\mapsto F_t}$, is a unitary representation of $\bf R$
on~$H$. We shall use Stone's theorem to produce a
(densely-defined unbounded) self-adjoint generator $F$. To
this end we must show that the representation, $t\mapsto F_t$,
is strongly continuous. In fact, we shall prove that it is
weakly continuous, which amounts to the same thing within the
set of unitaries.  For $a, b\in A$ and $t\in\bf R$, we have
\[(F_t \Lambda (a ), \Lambda (b))=(\Lambda (a), F_t^*\Lambda (b))
=(\Lambda (a), (\tilde{f}_{it}\otimes\iota )(W^*)\Lambda
(b))\]
\[=\tilde {f}_{it} ((\iota\otimes\omega_{\Lambda (a), \Lambda
(b)})(W^*)) =\tilde{f}_{it} \pi ((\iota\otimes\varphi
)((1\otimes b^* )\D (a)))\]
\[=f_{it}((\iota\otimes\varphi )((1\otimes b^* )\D (a))).\]
Hence, setting $c=(\iota\otimes\varphi )((1\otimes b^* )\D
(a))$, we get $(F_t \Lambda (a ), \Lambda (b)) =f_{it}(c)$.
The important point here is that $c$ is independent of~$t$.
Since $t\mapsto f_{it}(c)$ is continuous, the function
${t\mapsto (F_t \Lambda (a ), \Lambda (b))}$ is continuous.
That $t\mapsto (F_t u,v)$ is continuous, for all $u,v\in H$,
follows now by density of $\Lambda (A)$ in $H$.

That $(\tilde{f}_{it}\otimes\iota )\dar  (x) =F_t^* x F_t$,
for all $x\in \ar $, follows from the fact that $\dar  (x)=W^*
(1\otimes x)W$ and the extension ${\tilde{f}_{it}\otimes\iota
:M(\ar  \otimes B_0 (H))\rightarrow B(H)}$ is a
$*$-homomorphism and $1\otimes x\in  \ar   \otimes M(B_0
(H))\subseteq M(\ar \otimes B_0 (H))$. We showed above that
$(F_t \Lambda (a ), \Lambda (b))=f_{it}((\iota\otimes\varphi
)((1\otimes b^* )\D (a)))$. Taking $a=1$, it follows that
${(F_t \Lambda (1 ), \Lambda (b))}=f_{it}((\iota\otimes
\f)(1\otimes b^*))=\f(b^*)=\ip{\Lambda(1)}{\Lambda(b)}$, for
all $b\in A$. Hence, $F_t\Lambda(1)=\Lambda(1)$. Consequently,
$F_t\Lambda(a)=F_t\pi(a)\Lambda(1)=F_t\pi(a)F_t^*F_t\Lambda(1)=
(\tilde{f}_{-it}\otimes\iota )\dar
(\pi(a))\Lambda(1)=\pi(a*f_{-it})\Lambda(1)=\Lambda(a*f_{-it})$.

Let $\chi\colon A\otimes A\to A\otimes A$ be the flip map. Set
$\dop= \chi\D$. Then one easily verifies that $(A,\dop)$ is an
\aqg\ of compact type, called the {\em opposite} of $\AA$. The
Haar state of $(A,\dop)$ is the Haar state $\f$ of $\AA$.
Likewise, these quantum groups have the same co-unit. However,
the co-inverse of $(A,\dop)$ is $S^{-1}$, where $S$ is the
co-inverse of $\AA$. It is not hard to check that the family
of  modular functionals  associated to $(A,\dop)$ is the same
as that associated to $\AA$.

 Now let $\wop$ denote the
multiplicative unitary associated to $(A, \dop)$ in the GNS
construction $(H, \pi ,\Lambda )$ for the Haar state~$\f$. If
the co-multiplication on $\ar $ associated to $\dop$ on $A$ is
denoted by $\Delta_{\mathrm{op,\,r}}$, then it is easily
verified that $\Delta_{\mathrm{op,\,r}}(x)=\Sigma \dar
(x)\Sigma$, for all $x\in \ar $, where $\Sigma$ is the flip
operator on ${H\otimes H}$. 

 As before, define a
strongly-continuous one-parameter group $(E_t)_t$ of unitaries
by setting $E_t =(\tilde{f}_{it}\otimes\iota )(\wop)$, for all
$t\in\bf R$. Then $(\tilde{f}_{it}\otimes\iota
)\Delta_{\mathrm{op,\,r}} (x) =E_t^* x E_t$, for all $x\in
\ar $ and $E_t\Lambda(a)=\Lambda(f_{-it}*a)$, for all $a\in
A$. Also, there is a self-adjoint operator $E$ in $H$ such
that $E_t=\exp(i t E )$, for all~$t$. Clearly then
$(\iota\otimes \tilde{f}_{it})\dar  (x)=
(\tilde{f}_{it}\otimes\iota )\Delta_{\mathrm{op,\,r}}(x)
=E_t^* xE_t $.

Let $a\in A$. Then $\sigma_t (\pi (a))=\pi (f_{-it}*a*f_{-it}
)= {(\tilde{f}_{-it}\otimes \iota)\dar \pi(f_{-it}*a)}=F_t\pi
(f_{-it}*a )F^*_t =F_t(\iota\otimes \tilde{f}_{-it})\dar
\pi(a)F^*_t=F_t E_t\pi (a)E^*_t F^*_t $. It follows from
density of $\pi(A)$ in $\ar $ that $\s_t(x) ={F_t E_t xE^*_t
F^*_t}$, for all $x\in \ar $.

We defined the unitary $V_t$ on $H$ by setting $V_t\Lambda(a)=
\Lambda(f_{-it}*a*f_{-it})$, for all $a\in A$. Clearly, then,
$V_t=F_tE_t=E_tF_t$.

Recall that if $a\in A$, then $z\mapsto f_z(a)$ is an analytic
function on $\bC$. Set $\s_z(\pi(a))=\pi(f_{-iz}*a*f_{-iz})$.
Then it is easily verified that the map ${z\mapsto
\s_z\pi(a)}$ is analytic, in the sense that if $\t\in
A^*_{\mathrm{r}}$, then ${z\mapsto \t\s_z\pi(a)}$ is analytic.
Hence, the map ${z\mapsto \s_z\pi(a)}$ provides an analytic
extension to the plane of the function ${t\mapsto \s_t\pi(a)}$
on~$\bR$. This shows that $\pi(A)$ is contained in the set of
analytic elements for the C*-dynamical system $(\ar ,\s)$.
Moreover, it follows from Condition~5 in the list of
properties associated to the family $(f_z)_z$ that we stated
at the beginning of this section that $\f_{\mathrm{r}}(\pi(a)
\pi(b))=\f_{\mathrm{r}}(\pi(b)\s_i\pi(a))$. Now, as 
$\pi(A)$ is dense in $\ar$ and invariant under each $\sigma_{t}$ for 
$t \in \bR,$
$\pi(A)$ is a core for $\sigma_{i}$ (using \cite[Corollary 1.22]{Ku2}).
Hence, it
follows that  the state $\f_{\mathrm{r}}$ satisfies the KMS
condition for the automorphism group $(\s_t)_t$ at inverse
temperature $\beta=1$.

We summarize some of the previous discussion in the following
result.

\begin{thm}
For all $t\in\bf R$, set $F_t={(\tilde{f}_{it}\otimes\iota)(W)}$
and $E_t={(\tilde{f}_{it}\otimes\iota )(\wop)}.$ Then
$F_t\Lambda(a)=\Lambda(a*f_{-it})$ and
$E_t\Lambda(a)=\Lambda(f_{-it}*a)$, for all $a\in A$, from
which it follows that $F_tE_t=E_tF_t$, for all $t\in \bR$.
Moreover, there exists self-adjoint operators $F$ and $E$ in
$H$ such that $F_t=\exp(i t F )$ and $E_t=\exp(i t E )$, for
all $t\in \bR$.

For all $t\in \bR$, set  $V_t =F_tE_t$, so that $V_t=\exp(i t
(E+F))$. Then $\s_t(x)=V_tx V^*_t$, for all $x\in \ar $. Here
$\s_t$ is the unique automorphism of $\ar $ for which
$\s_t(\pi(a))=\pi (f_{-it}*a*f_{-it} )$, for all $a\in A$. The
Haar state $\f_{\mathrm{r}}$ on $\ar $ satisfies the KMS
condition for $(\s_t)_t$ at inverse temperature $\beta=1$.

\end{thm}

\bigskip
\noindent \emph{Acknowledgements}: The authors would like to thank the 
referee for his careful reading of the first version of this paper 
and for his many helpful comments.

\bigskip
{\parindent=0pt \parskip=4pt Addresses of the authors:

\smallskip Erik B\'edos, Institute of Mathematics, University of
Oslo, \\ P.B. 1053 Blindern, 0316 Oslo, Norway. E-mail:
bedos@math.uio.no.

Gerard J. Murphy, Department of Mathematics, National
University of Ireland, Cork, Ireland. E-mail: gjm@ucc.ie.

Lars Tuset, Faculty of Engineering, Oslo University College,
\\ Cort Adelers Gate~30, 0254 Oslo, Norway.
E-mail:Lars.Tuset@iu.hio.no.}

\end{document}